\documentclass{article}
\usepackage{amssymb}
\usepackage{amsmath,color}
\usepackage[hidelinks]{hyperref}

\usepackage[square,numbers]{natbib}         
\newtheorem{lem}{Lemma}[section]
\newtheorem{cor}[lem]{Corollary}
\newtheorem{teo}[lem]{Theorem}
\newtheorem{os}[lem]{Remark}

\newtheorem{prop}[lem]{Proposition}

\newcommand\sign{\mathop{\rm sign}}
\newcommand{\qed}{\thinspace\null\nobreak\hfill\hbox{\vbox{\kern-.2pt\hrule
 height.2pt depth.2pt\kern-.2pt\kern-.2pt \hbox to2.5mm{\kern-.2pt\vrule
 width.4pt \kern-.2pt\raise2.5mm\vbox to.2pt{}\lower0pt\vtop
 to.2pt{}\hfil\kern-.2pt \vrule
 width.4pt \kern-.2pt}\kern-.2pt\kern-.2pt\hrule height.2pt depth.2pt
 \kern-.2pt}}\par\medbreak}

\newcommand{\R}{\mathbb{R}}

\newcommand{\C}{\mathbb{C}}

\newcommand{\N}{\mathbb{N}}

\newcommand{\Rp}{\textrm{\emph{Re}\,}}
\newcommand{\Ip}{\textrm{\emph{Im}\,}}
\newcommand{\eps}{\varepsilon}

\newcommand{\ds}{\displaystyle}

\date{}
\begin{document}

\title{Singular parabolic problems in the half-space}
\author{G. Metafune \thanks{Dipartimento di Matematica e Fisica ``Ennio De Giorgi'', Universit\`a del Salento, C.P.193, 73100, Lecce, Italy.
-mail:  giorgio.metafune@unisalento.it}\qquad L. Negro \thanks{Dipartimento di Matematica e Fisica  ``Ennio De
Giorgi'', Universit\`a del Salento, C.P.193, 73100, Lecce, Italy e Dipartimento di Matematica, Universit\`a degli Studi di Salerno,  Via Giovanni Paolo II, 132 - 84084 Fisciano (SA). email: luigi.negro@unisalento.it} \qquad C. Spina \thanks{Dipartimento di Matematica e Fisica``Ennio De Giorgi'', Universit\`a del Salento, C.P.193, 73100, Lecce, Italy.
e-mail:  chiara.spina@unisalento.it}}

\maketitle
\begin{abstract}
\noindent 

We study elliptic and parabolic problems governed by singular elliptic   operators 
 \begin{equation*}
\mathcal L =\sum_{i,j=1}^{N+1}q_{ij}D_{ij}+\frac c y D_y
\end{equation*}
in the half-space $\R^{N+1}_+=\{(x,y): x \in \R^N, y>0\}$ under Neumann boundary conditions at $y=0$.  More general operators and oblique derivative boundary conditions are also considered.

\par

\end{abstract}

\section{Introduction}

 In this paper   we study  solvability and regularity of elliptic and parabolic problems associated to the  degenerate   operators 
 \begin{equation}\label{def L}
	\mathcal L =\sum_{i,j=1}^{N+1}q_{ij}D_{ij}+\frac c y D_y \quad {\rm and\ }\quad  D_t-\mathcal L
\end{equation}
in the half-space $\R^{N+1}_+=\{(x,y): x \in \R^N, y>0\}$ under Neumann boundary conditions at $y=0$.  Here $c\in\R$ and  $Q=\left(q_{ij}\right)_{i,j=1,\dots,N+1}$ is a constant, real, positive definite  matrix.

The special case where  
\begin{equation} \label{Lpart}
\mathcal L=\Delta_x+B_y, \qquad  B_y=D_{yy}+ \frac cy D_y
\end{equation}
 has been extensively studied in \cite{MNS-Caffarelli}. In this situation
$B_y$  is a one-dimensional Bessel operator and $\Delta_x$ and $B_y$ commute. These operators  play a major role in the investigation of the fractional powers $(-\Delta_x)^s$ and  $(D_t-\Delta_x)^s$, $s=(1-c)/2$, through the  ``extension procedure" of Caffarelli and Silvestre, see \cite{Caffarelli-Silvestre}. Nowadays this method has been extended to more general situations first in \cite{Stinga-Torrea-Extension} for self-adjoint operators, then in \cite{ GaleMiana-Extension} for generators of semigroup and in \cite{Arendt-ExtensionProblem} with precise regularity results  in the Hilbertian case.

The operator $\mathcal L$ is a model case of a more general situation which we explain below and will be treated in a  subsequent paper.

Let $\Omega$ be a smooth open bounded set in $\R^{k}$ and  let $d:=\mbox{dist}(\ \cdot\ ,\partial\Omega)$ be the distance from $\partial\Omega$. Let us consider operators of the form
\begin{equation*}      
	\mathcal L=\mbox{tr}\left(A(z)D^2\right)+ {\mbox{dist}(z,\partial\Omega)^{-1}} \left(b(z),\nabla\right)
	\qquad z \in \Omega,
\end{equation*}
under Neumann boundary condition at $\partial \Omega$. Here $A=\left(a_{ij}\right)$ is uniformly elliptic and $A, b\in C\left(\overline\Omega\right)$. Standard localization and freezing the coefficients reduce  the latter operator  to the half-space $\R^{N+1}_+$ in the form \eqref{def L}.  However, the form \eqref{Lpart} is not sufficient, since mixed derivatives and oblique boundary conditions appear in  the localization procedure, unless heavy restrictions are assumed.


\smallskip

Unfortunately dealing with such additional terms  is an hard complication of the problem, 
one important reason being the loss of commutativity. Indeed, generation properties and kernel estimates are easily proved for $\Delta_x+B_y$ by the same 
properties of the commuting blocks $\Delta_x$ and $B_y$, and this strategy clearly fails for operators like \eqref{def L}.

\smallskip
Moreover, the general  case cannot be reduced to the form \eqref{Lpart} by linear change of variables, as it is usually done for second order equations with constant coefficients. In fact, a linear change of variables which transforms $\sum_{ij}q_{ij}D_{ij}$ into $\Delta_x +D_{yy}$ acts also on $\frac cy D_y$ introducing an additional singular term like $\frac a y\cdot \nabla_x$, $ a \in \R^N$, and changing the Neumann boundary condition into an oblique one, see Section 8.

 The operators $\mathcal L$ in the general form \eqref{def L} have been considered in  \cite{dong2020parabolic} and \cite{dong2021weighted} and, in particular  in  \cite{dong2020RMI} with $|c|<1$ and in \cite{dong2020neumann} for $c>-1$, as in our assumptions. The authors show solvability and regularity of related elliptic and parabolic problems in weighted $L^p$-spaces, even for variable coefficients, using tools form linear PDE and Muckhenoupt weights.

In this paper we  use semigroup theory and operator valued harmonic analysis to prove similar results in  weighted $L^p$ spaces,  where the weight is a power  $y^m$, $m\in\R$. Even though the most important case are $m=0$, $m=c/\gamma$ where $\gamma=q_{N+1,N+1}$, which correspond to the Lebesgue and the symmetrizing measure, our methods work for all $m$ satisfying  $0<(m+1)/p<c+1$. 

In the language of semigroup theory, we prove that $\mathcal L$, endowed with  Neumann boundary conditions, generates an analytic semigroup in $L^p_m=L^p(\R^{N+1}_+; y^m dxdy)$,  we characterize its domain as a weighted Sobolev space and show that it has maximal regularity, which means that both $D_t v$ and $\mathcal L v$ have the same regularity as $(D_t -\mathcal L) v$. In comparison to \cite{dong2020neumann}, we obtain solvability of the problem $\lambda u-\mathcal L u=0$ also for complex $\lambda$.


 We prove both elliptic and parabolic estimates 
\begin{equation} \label{closedness}
\|D_{x_ix_j} u\|_{L^p_m}+\|D_{x_iy}u\|_{L^p_m}+\|D_{yy}u\|_{L^p_m}+\| y^{-1}D_y u\|_{L^p_m} \le C\| \mathcal L u\|_{L^p_m}, \quad 
\end{equation}
and
\begin{equation} \label{closedness1}
\|D_t u\|_{L^p_m}+\|\mathcal Lu\|_{L^p_m} \le C\| (D_t-\mathcal L) u\|_{L^p_m}, \quad 
\end{equation}
where the $L^p$ norms are taken over $\R_+^{N+1}$ and on $(0, \infty) \times \R_+^{N+1}$ respectively. 

%

In order to obtain \eqref{closedness} and  \eqref{closedness1} and show solvability, we use tools from vector-valued harmonic analysis: let us explain the main ideas. 

Assume, without loss of generality as explained later,  that $$\mathcal Lu=\Delta_x u+2a\cdot\nabla_xD_y u+B_yu=f.$$ Taking the Fourier transform with respect to $x$ we obtain $$-|\xi|^2 \hat u(\xi,y)+2ia\cdot\xi D_y+ B_y \hat u(\xi,y)=\hat f(\xi,y)$$ and then $$|\xi|^2 \hat u(\xi,y)=- |\xi|^2 (|\xi|^2-2ia\cdot\xi D_y-B_y)^{-1}\hat f (\xi,y), \quad\xi D_y\hat u(\xi,y)=\xi (|\xi|^2-2ia\cdot\xi D_y-B_y)^{-1}\hat f (\xi,y).$$
Denoting by ${\cal F}$ the Fourier transform with respect to $x$  we get 
$$
\Delta_x  \mathcal  L^{-1}=-{\cal F}^{-1}\left (|\xi|^2(|\xi|^2-2ia\cdot\xi D_y-B_y)^{-1} \right) {\cal F}
$$
and
$$
\nabla_x D_y  \mathcal  L^{-1}={\cal F}^{-1}\left (i\xi(|\xi|^2-2ia\cdot\xi D_y-B_y)^{-1} \right) {\cal F}.
$$
The boundedness of $\Delta_x  \mathcal L^{-1}$ is equivalent to that of the multiplier $\xi \in \R^N \to |\xi|^2(|\xi|^2-2ia\cdot\xi D_y-B_y)^{-1} $ in $L^p(\R^N; L^p_m(0,\infty))=L^p_m$.  Similarly,  the boundedness of $\nabla_xD_y  \mathcal L^{-1}$ is equivalent to that of the multiplier $\xi \in \R^N \to \xi(|\xi|^2-2ia\cdot\xi D_y-B_y)^{-1} $.

To prove solvability one more multiplier is needed, namely $\xi \in \R^N \to \lambda(\lambda+|\xi|^2-2ia\cdot\xi D_y - B_y)^{-1}$, where  $\lambda$ belongs to a sector in the complex plane.

The proof of the  boundedness of these multipliers uses a vector valued Mikhlin multiplier theorem which rests on square function estimates. 
The strategy for proving  \eqref{closedness1} is similar after taking the Fourier transform with respect to $t$.

We refer to  \cite{KunstWeis},  \cite{DenkHieberPruss} and the new books \cite{WeisBook1}, \cite{WeisBook2} for  this approach, which we recall in Section 2. 

A first crucial step consists in the study of the one-dimensional operator $L=B_y+ib D_y$, which is of independent interest. 
Since $b=2a \cdot \xi$  in the multipliers above, we need precise dependence on $b$ in all estimates.

%

\medskip
The paper is organized as  follows. In Section 2 we briefly recall the tools of vector-valued harmonic we need.

In Sections 3 and  4 we recall some results concerning  weighted Sobolev spaces and the  one-dimensional Bessel operator $B_y$.

In Section 5 we define the $1d$-operator $L$  through a quadratic form in  $L^2((0,+\infty), y^c\, dy)$. We prove heat kernel estimates for real times by domination  and then we extend them to complex times, via Davies Gaffney estimates.

In Section 6 we prove the boundedness of the multipliers introduced above.

In Section 7, which is the core of the paper,  we prove generation results, maximal regularity and domain characterization for the operator $\mathcal L$, under Neumann boundary conditions. 

Finally, in Section 8, we extend our results by considering operators of the form 
\begin{equation*}
\mathcal L =\sum_{i,j=1}^{N+1}q_{ij}D_{ij}+\frac c y D_y+\frac{a\cdot\nabla_x}{y},
\end{equation*}
and oblique derivative boundary conditions.

\bigskip
\noindent\textbf{Notation.} For $N \ge 0$, $\R^{N+1}_+=\{(x,y): x \in \R^N, y>0\}$. We write $\nabla u, D^2 u$ for the gradient and the Hessian matrix of a function $u$ with respect to all $x,y$ variables and $\nabla_x u, D_y u, D_{x_ix_j }u, D_{x_i y} u$ and so on, to distinguish the role of $x$ and $y$.

For $m \in \R$ we consider the measure $y^m dx dy $ in $\R^{N+1}_+$ and  we write $L^p_m(\R_+^{N+1})$, and often only $L^p_m$ when $\R^{N+1}_+$ is understood, for  $L^p(\R_+^{N+1}; y^m dx dy)$. 

Similarly $W^{k,p}_m=W^{k,p}_m(\R^{N+1}_+)=\{u \in L^p_m: \partial^\alpha u \in  L^p_m,  \quad |\alpha| \le k\}$. We use $\hat {C}^k([0,\infty[)$ for the space of uniformly continuous, $k$-times differentiable functions on $[0,\infty[$, tending to zero at infinity with all derivatives. 

$\C^+$ stands for $\{ \lambda \in \C: \Rp \lambda >0 \}$ and, for $|\theta| \leq \pi$, we denote by  $\Sigma_{\theta}$  the open sector $\{\lambda \in \C: \lambda \neq 0, \ |Arg (\lambda)| <\theta\}$.  

Given $a$ and $b$ $\in\R$, $a\wedge b$, $a \vee b$  denote  their minimum and  maximum. We  write $f(x)\simeq g(x)$ for $x$ in a set $I$ and positive $f,g$, if for some $C_1,C_2>0$ 
\begin{equation*}
	C_1\,g(x)\leq f(x)\leq C_2\, g(x),\quad x\in I.
\end{equation*}

\bigskip

\section{Vector-valued harmonic analysis}
Regularity properties for $\mathcal L$ and $D_t-\mathcal L$ follow once we prove the estimate
\begin{equation}  \label{regularity}
	\|D^2 u\|_p \le C \| \mathcal L u\|_p, \quad \|D_t u\|_p+\|D^2 u\|_p \le C \|(D_t- \mathcal L) u\|_p
\end{equation}
for $u$ in an appropriate Sobolev space and this is equivalent to saying that $D^2 {\mathcal L}^{-1}$ and $D_t(D_t-\mathcal L)^{-1} $ are bounded operators. We use  the strategy which arose first in the study of maximal regularity of parabolic problems, that is for the equation $u_t=Au+f, u(0)=0$ where $A$ is the generator of an analytic semigroup on a Banach space $X$. Estimates like
$$
\|u_t\|_p+\|Au\|_p \le \|f\|_p
$$
arr equivalent to the  boundedness of the operator $A(D_t-A)^{-1}$.

This strategy  relies on Mikhlin vector-valued multiplier theorems which we now recall here,
 referring the reader to \cite{DenkHieberPruss}, \cite{Pruss-Simonett} or \cite{KunstWeis} for all proofs.

Let ${\cal S}$ be a subset of $B(X)$, the space of all bounded linear operators on a Banach space $X$. ${\cal S}$ is $\mathcal R$-bounded if there is a constant $C$ such that
$$
\|\sum_i \eps_i S_i x_i\|_{L^p(\Omega; X)} \le C\|\sum_i \eps_i  x_i\|_{L^p(\Omega; X)} 
$$
for every finite sum as above, where $(x_i ) \subset X, (S_i) \subset {\cal S}$ and $\eps_i:\Omega \to \{-1,1\}$ are independent and symmetric random variables on a probability space $\Omega$. In particular ${\cal S}$ is a bounded subset of $B(X)$. The smallest constant $C$ for which the above definition holds is the $\mathcal R$-bound of $\mathcal S$,  denoted by $\mathcal R(\mathcal S)$.
It is well-known that this definition does not  depend on $1 \le p<\infty$ (however, the constant $\mathcal R(\mathcal S)$ does) and that $\mathcal R$-boundedness  is equivalent to boundedness when $X$ is an Hilbert space.
When $X$ is an $L^p(\Sigma)$ space (with respect to any $\sigma$-finite measure defined on a $\sigma$-algebra $\Sigma$), testing  $\mathcal R$-boundedness is equivalent to proving square functions estimates, see \cite[Remark 2.9]{KunstWeis}.

\begin{prop}\label{Square funct R-bound} Let ${\cal S} \subset B(L^p(\Sigma))$, $1<p<\infty$. Then ${\cal S}$ is $\mathcal R$-bounded if and only if there is a constant $C>0$ such that for every finite family $(f_i)\in L^p(\Sigma), (S_i) \in {\cal S}$
	$$
	\left\|\left (\sum_i |S_if_i|^2\right )^{\frac{1}{2}}\right\|_{L^p(\Sigma)} \le C\left\|\left (\sum_i |f_i|^2\right)^{\frac{1}{2}}\right\|_{L^p(\Sigma)}.
	$$
\end{prop}
The best constant $C$ for which the above square functions estimates hold satisfies $\kappa^{-1} C \leq \mathcal R(\mathcal S) \leq \kappa C$ for a suitable $\kappa>0$ (depending only on $p$).  Using the proposition above, $\mathcal R$-boundedness follows from domination by a positive $\mathcal R$-bounded family.
\begin{cor} \label{domination}
	 Let  ${\cal S}, {\cal T} \subset B(L^p(\Sigma))$, $1<p<\infty$ and assume that $\cal T$ is an $\mathcal R$ bounded family of positive operators and that for every $S \in \cal S$ there exists $T \in \cal T$ such that $|Sf| \leq T|f|$ pointwise, for every $f \in L^p(\Sigma)$. Then ${\cal S}$ is $\mathcal R$-bounded.
\end{cor}

We also need the following result about the integral mean of a $\mathcal R$-bounded  family of operator which we state in the version we use.
\begin{prop}{\em \cite[Corollary 2.14]{KunstWeis}} \label{Mean R-bound}
	Let  $X$ be a Banach space and let  ${\cal F}\subset B(X)$ be an $\mathcal R$-bounded  family of operator. For every strongly measurable $N:\Sigma\to B(X)$ on a $\sigma$-finite measure space  $(\Sigma,\mu)$ with values in 
	$\cal F$ and every $h\in L^1\left(\Sigma,\mu \right)$ we define the operator  $T_{N,\cal F}\in B(X)$ by
	\begin{align*}
	T_{N,\cal F} x=\int_{\Sigma} h(\omega)N(\omega)xd\mu(\omega),\qquad x\in X.
	\end{align*}
Then the family 
\begin{align*}
	\mathcal C=\left\{T_{N,\cal F}\;:\; \|h\|_{L^1}\leq 1, N\;\text{as above}\right\}
\end{align*} is $\mathcal R$ bounded and $\mathcal R(\mathcal C)\leq 2\mathcal R(\mathcal F)$.
\end{prop}
Let $(A, D(A))$ be a sectorial operator in a Banach space $X$; this means that $\rho (-A) \supset \Sigma_{\pi-\phi}$ for some $\phi <\pi$ and that $\lambda (\lambda+A)^{-1}$ is bounded in $\Sigma_{\pi-\phi}$. The infimum of all such $\phi$ is called the spectral angle of $A$ and denoted by $\phi_A$. Note that $-A$ generates an analytic semigroup if and only if $\phi_A<\pi/2$. The definition of $\mathcal  R$-sectorial operator is similar, substituting boundedness of $\lambda(\lambda+A)^{-1}$ with $\mathcal R$-boundedness in $\Sigma_{\pi-\phi}$. As above one denotes by $\phi^R_A$ the infimum of all $\phi$ for which this happens; since $\mathcal R$-boundedness implies boundedness, we have $\phi_A \le \phi^R_A$.

\medskip

The $\mathcal R$-boundedness of the resolvent characterizes the regularity of the associated inhomogeneous parabolic problem, as we explain now.

An analytic semigroup $(e^{-tA})_{t \ge0}$ on a Banach space $X$ with generator $-A$ has
{\it maximal regularity of type $L^q$} ($1<q<\infty$)
if for each $f\in L^q([0,T];X)$ the function
$t\mapsto u(t)=\int_0^te^{-(t-s)A})f(s)\,ds$ belongs to
$W^{1,q}([0,T];X)\cap L^q([0,T];D(A))$.
This means that the mild solution of the evolution equation
$$u'(t)+Au(t)=f(t), \quad t>0, \qquad u(0)=0,$$
is in fact a strong solution and has the best regularity one can expect.
It is known that this property does not depend on $1<q<\infty$ and $T>0$.
A characterization of maximal regularity is available in UMD Banach spaces, through the $\mathcal  R$-boundedness of the resolvent in a suitable sector $\omega+\Sigma_{\phi}$, with $\omega \in \R$ and $\phi>\pi/2$ or, equivalently, of the scaled semigroup $e^{-(A+\omega')t}$ in a sector around the positive axis. In the case of $L^p$ spaces it can be restated in the following form,  see \cite[Theorem 1.11]{KunstWeis}

\begin{teo}\label{MR} Let $(e^{-tA})_{t \ge0}$ be a bounded analytic semigroup in $L^p(\Sigma)$, $1<p<\infty$,  with generator $-A$. Then $T(\cdot)$ has maximal regularity of type $L^q$  if and only if the set $\{\lambda(\lambda+A)^{-1}, \lambda \in  \Sigma_{\pi/2+\phi} \}$ is $\mathcal R$- bounded for some $\phi>0$. In an equivalent way, if and only if 
	there are constants $0<\phi<\pi/2 $, $C>0$ such that for every finite sequence $(\lambda_i) \subset \Sigma_{\pi/2+\phi}$, $(f_i) \subset  L^p$
	$$
	\left\|\left (\sum_i |\lambda_i (\lambda_i+A)^{-1}f_i|^2\right )^{\frac{1}{2}}\right\|_{L^p(\Sigma)} \le C\left\|\left (\sum_i |f_i|^2\right)^{\frac{1}{2}}\right\|_{L^p(\Sigma)}
	$$
	or, equivalently, there are constants $0<\phi'<\pi/2 $, $C'>0$ such that  for every finite sequence
	$(z_i) \subset \Sigma_{\phi'}$, $(f_i) \subset  L^p$
	$$
	\left\|\left (\sum_i |e^{-z_i A}f_i|^2\right )^{\frac{1}{2}}\right\|_{L^p(\Sigma)} \le C'\left\|\left (\sum_i |f_i|^2\right)^{\frac{1}{2}}\right\|_{L^p(\Sigma)}.
	$$
\end{teo}
\medskip

Finally we state  a version of the operator-valued Mikhlin multiplier theorem in the N-dimensional case, see e.g. \cite[Corollary 8.3.22]{WeisBook2}.

\begin{teo}   \label{marcinkiewicz}
	Let $1<p<\infty$, $M\in C^N(\R^N\setminus \{0\}; B(L^p(\Sigma))$ be such that  the set
	$$\left \{\xi^{\alpha}D^\alpha_\xi M(\xi): \xi\in \R^{N}\setminus\{0\}, \ \alpha\in\{0,1\}^N \right \}$$
	is $\mathcal{R}$-bounded.
	Then the operator $T_M={\cal F}^{-1}M {\cal F}$ is bounded in $L^p(\R^N, L^p(\Sigma))$, where $\cal{F}$ denotes the Fourier transform.
\end{teo}
\section{Weighted Sobolev spaces}

We collect in this section  the main results concerning  weighted Sobolev spaces, referring to \cite{MNS-Sobolev} for further details  and proofs, even in more general situations. 

For $p>1$, $m >-1$ , we define the Sobolev space
\begin{align*}
	W^{2,p}_m &=\left\{u\in W^{2,p}_{loc}(\R^{N+1}_+):\ u,\  D_{x_i}u, D_y u,  D_{x_ix_j}u,\  D_{x_i y}u,  D_{yy}u \in L^p_m\right\}
\end{align*}
which is a Banach space equipped with the norm
\begin{align*}
	\|u\|_{W^{2,p}_m}=&\|u\|_{L^p_m}+\sum_{i=1}^N\| D_{x_i}u\|_{L^p_m}+\|D_{y}u\|_{L^p_m}+\sum_{i,j=1}^N\| D_{x_ix_j}u\|_{L^p_m}+\|D_{yy}u\|_{L^p_m}+\sum_{i=1}^N\| D_{x_iy}u\|_{L^p_m}.
\end{align*}
Next we add a Neumann boundary condition for $y=0$  in the form $y^{-1}D_yu\in L^p_m$ and set
\begin{align*}
	W^{2,p}_{m, \mathcal N}=\{u \in W^{2,p}_m:\  y^{-1}D_yu\ \in L^p_m\}
\end{align*}
with the norm
$$
\|u\|_{W^{2,p}_{m,\mathcal N}}=\|u\|_{W^{2,p}_m}+\|y^{-1}D_yu\|_{ L^p_m}.
$$

\begin{os}\label{Os Sob 1-d}
	With obvious changes we consider also the analogous Sobolev spaces $W^{2,p}_m$ and $W^{2,p}_{m,\cal N}$ on $\R_+$. 
\end{os}

The next result clarifies in which sense the condition $y^{-1}D_y u \in L^p_m$ is a Neumann boundary condition.

\begin{prop} \label{neumann} The following assertions hold.
	\begin{itemize} 
		\item[(i)] If $\frac{m+1}{p} >1$, then $W^{2,p}_{m,\mathcal N}=W^{2,p}_m$.
		\item[(ii)] If $\frac{m+1}{p} <1$, then $$W^{2,p}_{m,\mathcal N}=\{u \in W^{2,p}_m: \lim_{y \to 0}D_yu(x,y)=0\ {\rm for\ a.e.\   x \in \R^N }\}.$$
	\end{itemize}
	In both cases (i) and (ii), the norm of $W^{2,p}_{m,\mathcal N}$ is equivalent to that of $W^{2,p}_m$.
\end{prop}

We provide  an equivalent description of $W^{2,p}_{m,\mathcal N}$, adapted to the operator $D_{yy}+cy^{-1}D_y$.
\begin{prop}\label{Trace D_yu in W}
	Let   $0<\frac{m+1}{p}<c+1$.  Then
	\begin{align*} 
		W^{2,p}_{m,\mathcal N}=&\left\{u \in  W^{2,p}_{loc}(\R^{N+1}_+): u,\ \Delta_xu\in L^p_m,  D_{yy}u+c\frac{D_yu}y \in L^p_m  \right. \\[1ex]			
		&\left.\hspace{10ex}\text{\;\;and\;\;}\lim_{y\to 0}y^c D_yu(x,y)=0\ {\rm for\ a.e.\   x \in \R^N }\}\right\}
	\end{align*}
	and the norms $\|u\|_{W^{2,p}_{m,\mathcal N}}$ and $$\|u\|_{L^p_m}+\|\Delta_x u\|_{L^p_m}+\|(D_{yy}u+cy^{-1}D_yu)\|_{L^p_m}$$ are equivalent on $W^{2,p}_{m,\mathcal N}$.
\end{prop}

The next results  show the density  of smooth functions in $W^{2,p}_{m,\mathcal N}$. Let
\begin{equation} \label {defDcore}
	\mathcal{D}=\left \{u \in C_c^\infty ([0, \infty)), \ D_y u(y)=0\  {\rm for} \ y \leq \delta\ {\rm  and \ some\ } \delta>0\right \}
\end{equation}
and  (finite sums below)
\begin{equation}\label{defDcore1}
C_c^\infty (\R^{N})\otimes\mathcal D=\left\{u(x,y)=\sum_i u_i(x)v_i(y), \  u_i \in C_c^\infty (\R^N), \  v_i \in \cal D \right \}.
\end{equation}
\begin{teo} \label{core gen}
	If  $\frac{m+1}{p}>0$
	then $C_c^\infty (\R^{N})\otimes\mathcal D$ is dense in $W^{2,p}_{m,\mathcal N}$.
\end{teo}

Note that the condition $m+1>0$  is necessary for the inclusion  $C_c^\infty (\R^{N})\otimes\mathcal D \subset W^{2,p}_{m,\mathcal N}$. 
\medskip



\section{The Bessel operator $B$}\label{sec y^aB}
In this section we consider for   $c>-1$  the  one-dimensional Bessel operator 
\begin{equation*} 
	 B =D_{yy}+\frac{c}{y}D_y 
\end{equation*}
in the space $L^p_m=L^p_m(\R_+)$, under Neumann boundary conditions at $y=0$. We summarize below all the main results we need,  referring the reader to  \cite[Section 3]{MNS-Caffarelli}, \cite[Section 4]{MNS-PerturbedBessel}, \cite{MNS-CompleteDegenerate}  for further details  and to  \cite{MNS-Max-Reg,MNS-Grad, MNS-Sharp,Negro-Spina-Asympt} for analogous results in the multidimensional case.

Setting $H^1_c=\{u \in L^2_c, u' \in L^2_c\}$,  the operator $B$ (with Neumann boundary conditions)  is associated to the non-negative, symmetric and closed form  in $L^2_c$
\begin{align*}
\mathfrak{a}(u,v)
&:=
\int_0^\infty D_y u  D_y \overline{v}\,y^c dy,\qquad  D(\mathfrak{a})=H^1_c.
\end{align*}

If $1<p<\infty$, we recall that 
\begin{equation} \label{w2p}
	W^{2,p}_{\mathcal N}(m)=\left\{u\in W^{2,p}_{loc}(\R_+):\ u,\    D_{yy}u,\ D_{y}u,\ y^{-1}D_{y}u\in L^p_m\right\}
\end{equation}
and for $p=\infty$ we define
\begin{equation} \label{c2n}
	\hat{C}^2_{\mathcal N}=\left\{u\in \hat C^2([0, \infty[): D_yu(0)=0 \right\}
\end{equation}
(we recall that the "hat"  means that the functions and the derivatives tend to 0 at infinity, see the Notation). 
The Neumann boundary condition, denoted by the pedix $\mathcal N$,  is enclosed in the requirement $y^{-1}D_{y}u\in L^p_m$. This last is redundant when $(m+1)/p>1$ and equivalent to $D_yu(y) \to 0$ as $y \to 0$, when $(m+1)/p<1$, see Proposition  \ref{neumann}.

Consequently we write $ B^n$ or, more pedantically $ B^n_{m,p}$ if necessary, for the operator $ B$ endowed with the domain $W^{2,p}_{\mathcal N}(m)$.
This time the suffix $n$ reminds the Neumann boundary condition at $y=0$.

\begin{teo}\label{Neumnan alpha m 2}
If  $1<p<\infty$ and $0<\frac{m+1}p<c+1$, then 
$ B$ endowed with domain $W^{2,p}_{\mathcal N}(m)$
generates a bounded analytic semigroup $e^{zB}$ of angle $\pi/2$ on $L^p_m$, which is positive for $z>0$.

If $p =\infty$, then $B$ with domain $\hat{C}^2_{\mathcal N}$ generates an analytic semigroup of angle $\frac \pi 2$ in $\hat {C}([0, \infty[)$.
\end{teo}
{\sc Proof.} See \cite[Proposition 3.3, Proposition 5.3]{MNS-Caffarelli} and \cite[Proposition 3.4]{MNS-Caf-Schauder} for $p=\infty$.
\qed



\begin{prop}\label{Estimates Bessel kernels}
Let $c>-1$. The semigroup $(e^{zB^n})_{z \in \C_+}$ consists of integral operators.  Its heat kernel $p_B$, written  with respect the measure $\rho^cd\rho$, satisfies for every $\eps>0$, $z\in\Sigma_{\frac{\pi}{2}-\eps}$ and some  $C_\eps, \kappa_\eps>0$, 
	\begin{align*} 
		| p_{B}(z,y,\rho)|
		&\leq 
		C_\eps |z|^{-\frac{1}{2}}\rho^{-c}  \left (\frac{\rho}{|z|^{\frac{1}{2}}}\wedge 1 \right)^{c}
		\exp\left(-\frac{|y-\rho|^2}{\kappa_\eps |z|}\right),\\[1ex] 
		| D_y p_{B}(z,y,\rho)|
		&\leq 
		C_\eps |z|^{-1}\rho^{-c}\left (\frac{y}{|z|^{\frac{1}{2}}}\wedge 1 \right)\left (\frac{\rho}{|z|^{\frac{1}{2}}}\wedge 1 \right)^{c}
		\exp\left(-\frac{|y-\rho|^2}{\kappa_\eps |z|}\right).
	\end{align*}
\end{prop}
{\sc{Proof.}} See \cite[Propositions 2.8 and 2.9]{MNS-Caffarelli}. \qed 
\subsection{$\mathcal{R}$-boundedness of some families of operators associated to $B^n$}

 We consider a two-parameters  family of integral operators $\left(S_{\alpha,\beta}(t)\right)_{t>0}$ on $L^p_m$, defined for $\alpha,\beta\in\R$  and $t>0$ by
\begin{align*}
	S^{\alpha,\beta}(t)f(y)=t^{-\frac {1} 2}\,\left (\frac{y}{\sqrt t}\wedge 1 \right)^{-\alpha}\int_{0}^\infty   \left (\frac{z}{\sqrt t}\wedge 1 \right)^{-\beta}
	\exp\left(-\frac{|y-z|^2}{\kappa t}\right)f(z) \,dz,
\end{align*}
where $\kappa$ is a positive constant. We omit the dependence on $\kappa$ even though in some proofs we need to vary it. 
We also define the families 
\begin{equation} \label{defgammabeta}
\Gamma (\lambda)=\int_0^\infty  e^{-\lambda t}S^{0,-c}(t)\, dt, \qquad 
\Psi (\lambda)=\int_0^\infty  \frac{e^{-\lambda t}}{\sqrt{t}}S^{0,-c}(t)\, dt, \quad \lambda>0
\end{equation}
By Proposition \ref{Estimates Bessel kernels} and the results in 
 \cite{MNS-Caffarelli, MNS-Sharp},  $S^{0,-c}(t)$ and $\Gamma(\lambda)$ sharply estimates the behaviour of the semigroup $e^{tB^n}$ and of the resolvent of $B^n$,  in the sense that, for $f \geq 0$, 
\begin{align*}
e^{tB^n}f&\simeq S^{0,-c}(t)f,\quad t>0,\qquad (\lambda-B^n)^{-1}f\simeq \Gamma(\lambda)f,\quad \lambda>0.
\end{align*}
Similarly, $S^{0,-c}(t)/ \sqrt{t}$ and  $\Psi(\lambda)$   estimate the behaviour of the spatial derivative of the semigroup and of  the resolvent, that is 
\begin{align*}
|D_ye^{tB^n}f| \leq C \frac{S^{0,-c}(t)}{\sqrt t}|f|,\quad t>0,\qquad |D_y(\lambda-B^n)^{-1}f| \leq C \Psi(\lambda)|f|,\quad \lambda>0.
\end{align*}

We summarize in the following proposition the main properties of the above families, referring to \cite[Section 7]{MNS-Caffarelli} and \cite[Section 4]{MNS-Max-Reg} for further details.
\begin{prop}\label{R-bound S^a,b}
Let  $1<p< \infty$ and $\alpha,\,\beta\in\R$. The following properties hold.
\begin{itemize}
	\item[(i)] If $\alpha<\frac{m+1}p < 1-\beta$, then the family
	$\left(S^{\alpha,\beta}(t)\right)_{t\geq 0}$ is $\mathcal{R}$-bounded on $L^p_m$.
	\item[(ii)] If $0<\frac{m+1}p<c+1$ then the family
	$\left(\lambda\Gamma(\lambda)\right)_{\lambda>0}$ is $\mathcal{R}$-bounded on $L^p_m$.
	\item[(iii)] If $0<\frac{m+1}p<c+1$ then the family
	$\left(\sqrt\lambda\Psi(\lambda)\right)_{\lambda>0}$ is $\mathcal{R}$-bounded on $L^p_m$.
\end{itemize}
\end{prop}
{\sc{Proof.}} Property (i) follows from \cite[Theorem 7.7]{MNS-Caffarelli} and  (ii)  by (i) and  Proposition \ref{Mean R-bound}, with $h_\lambda(t)=\lambda e^{-\lambda t}$.
Property (iii) follows similarly by (i) and  Proposition \ref{Mean R-bound}, with $g_\lambda(t)=\sqrt{\frac{\lambda}{ t}} e^{-\lambda t}$.
\qed

\section{The operator $L_b=B+ibD_y-\frac{b^2}{4}$}\label{section L_b}
In this section we prove  generation properties and  heat kernel bounds for the operator   $$L_b:=B+ibD_y-\frac{b^2}{4}, $$ following the method of \cite[Sections 3, 4]{met-calv-negro-spina}. Note that $L_0=B$.

Recalling that  the operator $B^n$   is associated to the form  in $L^2_c$
\begin{align*}
\mathfrak{a}(u,v)
&:=
\int_0^\infty D_y u  D_y \overline{v}\,y^c dy,\qquad  D(\mathfrak{a})=H^1_c, 
\end{align*}
we  consider the perturbed  form $\mathfrak{ a}_b$  defined on $ D(\mathfrak{a}_b)=H^1_c$ by
\begin{equation}\label{form v}
 \mathfrak{a}_b(u,v)=\mathfrak{a}(u,v)- i b \left\langle D_yu,v\right\rangle_{L^2_c}+\frac{b^2}{4} \left\langle u,v\right\rangle_{L^2_c}
\end{equation}
and define $L_b$ in $L^2_c$ as the operator  associated to the form $\mathfrak{a}_b$, that is
\begin{align} \label{defLb}
\nonumber D(L_b)&=\{u \in D(\mathfrak{a}_b): \exists  f \in L^2_c \ {\rm such\ that}\  \mathfrak{a}_b(u,v)=\int_0^\infty f \overline{v}y^c\, dy\ {\rm for\ every}\ v\in D(\mathfrak{a}_b)\},\\[1ex]
  L_b u&=-f.
\end{align}

\subsection{The auxiliary operator $B-i\frac{bc}{2y}$}

For technical reasons we  consider also the form  $\mathfrak{\tilde a}_b$  defined on $D(\mathfrak{\tilde a}_b)=H^1_c$  by
\begin{equation}\label{form tilde}
	\nonumber\mathfrak{\tilde a}_b(u,v)=\mathfrak{a}(u,v)- i\frac  b 2\left\langle u,D_yv\right\rangle_{L^2_c}- i\frac  b 2\left\langle D_yu,v\right\rangle_{L^2_c}=
	\int_{\R_+}
	\left(D_y u D_y \overline{v}- i\frac  b 2D_y\left(u\overline{v}\right)\right)y^c\,dy
\end{equation}
and its   associated operator $A_b$ in $L^2_c$.
Since for smooth functions with compact support away from the origin
$$
c\int_0^\infty u \overline{v} y^{c-1}\, dy=
 \int_0^\infty u\big  (D_y (y^c\overline{v})-y^c  D_y \overline{v})\big )\, dy=- \int_0^\infty D_y(u\overline{v})y^c\, dy,
$$
the operator $A_b$ is defined on smooth functions by 
\begin{align*}
	A_b:= B-i\frac{bc}{2y}, \quad  A_0=B^n.
\end{align*}

We collect in the following proposition the main property satisfied by $\mathfrak{\tilde a}_b$.

\begin{prop}\label{prop form a_b tilde}
	If $c+1>0$,   then $\mathfrak{\tilde a}_b$ 
	is  accretive  and closed in $L^2_{c}$. Moreover
	\begin{itemize}
		\item[(i)] the adjoint form $\mathfrak{\tilde a}_b^\ast:(u,v)\mapsto \overline{\mathfrak{\tilde a}_b(v,u)}$ satisfies $\mathfrak{\tilde a}_b^\ast =\mathfrak{\tilde a}_{-b}$; 
		\item[(ii)] its real part  is $\Rp{\mathfrak{\tilde a}_b}:=\frac{\mathfrak{\tilde a}_b+\mathfrak{\tilde a}_b^\ast}{2}=\mathfrak{ a}$; 
		\item[(iii)] for any $\epsilon>0$, $u \in H^1_c$
		\begin{align*}
			\left|\Ip{\mathfrak{\tilde a}_b}(u,u)\right|
			\leq \eps \left(\mathfrak{ a}(u,u)+\frac{b^2}{4\epsilon^2}\|u\|^2_{L^2_c}\right ).
		\end{align*}
	\end{itemize}
\end{prop}
{\sc Proof.}  Properties (i) and (ii) are immediate consequences of the definition.	Since $\Rp{\mathfrak{\tilde a}_b(u,u)}=\mathfrak{ a}(u,u)\geq 0$,    $\mathfrak{\tilde a}_b$ 
is  accretive and, furthermore,  the norm induced by the form $\mathfrak{\tilde a}_b$  coincides with the one induced by $\mathfrak{ a}$ and then $\mathfrak{\tilde a}_b$ is closed.  

To prove   (iii)   it is enough to observe that for any $\epsilon>0$ and $u\in H^1_c$ one has 
\begin{align*}
	\left|\Ip{\mathfrak{\tilde a}_b}(u,u)\right|&=\left|-b\int_0^\infty \Rp\left(\overline u D_yu\right)y^c\,dy\right|\leq |b| \left(\|D_yu\|_{L^2_c} \|u\|_{L^2_c}\right)\\[1ex]
	&=\epsilon  \left(2\,\|D_yu\|_{L^2_c} \left|\frac{b}{2\epsilon}\right|\|u\|_{L^2_c}\right)\leq \epsilon \left(\mathfrak{ a}(u,u)+\frac{b^2}{4\epsilon^2}\|u\|^2_{L^2_c}\right)
\end{align*}
where we used the elementary inequality $D_y(|u|^2)=2\Rp\left(\overline u D_yu\right)$.
\qed
By standard theory on sesquilinear forms we have  the following results.

\begin{prop}\label{kernel Vreal tilde}
	If $c+1>0$,  then the   operator $A_b$ generates  an  analytic semigroup  of angle $\frac\pi 2$ in $L^2_c$   which satisfies for any $\epsilon>0$
	\begin{align*}
		\|e^{zA_b}f\|_{L^2_c}\leq e^{\frac{b^2}{4\epsilon^2}\Rp z}\|f\|_{L^2_c},\qquad \forall z\in \Sigma_{\frac\pi 2-\arctan \epsilon }.
	\end{align*} 
	Moreover
	\begin{itemize}
		\item[(i)] The semigroup $\left(e^{tA_b}\right)_{t\geq 0}$ is  $L^\infty$-contractive  and it is dominated by $e^{tB^n}$, that is 
		\begin{align*}
			|e^{tA_b}f|\leq e^{tB^n}|f|,\qquad t>0,\quad  f\in L^2_c.
		\end{align*}
		\item[(ii)] $\left(e^{tA_b}\right)_{t\geq 0}$ is a semigroup of integral operators and its heat kernel $\tilde p_b$, taken with respect to the measure $\rho^cd\rho$, satisfies for some constant $C$ independent of $b$
		\begin{align*}
			|\tilde p_{b}(t,y,\rho)|\leq C t^{-\frac{1}{2}}  \rho^{-c}\left (\frac{\rho}{t^{\frac{1}{2}}}\wedge 1 \right)^{c}
			\exp\left(-\frac{|y-\rho|^2}{\kappa t}\right),\qquad \text{for a.e. $y,\rho>0$}.
		\end{align*}     
		\item[(iii)] ${A}^\ast_b=\tilde A _{-b}$ and for any $s>0$ the operator satisfies the scaling property
		\begin{align*}
			I_{\frac 1 s}\circ A_b \circ I_s= s^2A_{\frac b s},\qquad I_s u(y):=u(sy). 
		\end{align*}
	\end{itemize}  
\end{prop}
{\sc Proof.} The generation properties  follows  using  Proposition \ref{prop form a_b tilde} and \cite[Teorems 1.52, 1.53]{Ouhabaz}.
$e^{tA_b}$ is $L^\infty$-contractive    from  \cite[Theorem 2.13]{Ouhabaz}.

The domination property follows  from \cite[Corollary 2.21]{Ouhabaz}.
(ii) is a consequence of  \cite[Proposition 1.9]{ArendtBukhalov} since 
$e^{tA_b}$ is dominated by the positive integral operator $e^{tB^n}$ whose kernel satisfies the stated estimate, see  \cite[Proposition 2.8]{MNS-Caffarelli} where, however, the  kernel is written with respect to the Lebesgue measure. (iii) follows from (i) of Proposition \ref{prop form a_b tilde} and by elementary computation.\qed

\subsection{Bounds for $e^{t L_b}$}

The following elementary lemma relates $L_b$ to $A_b$.
\begin{lem} \label{iso}
The isometry 
\begin{align*}
	T:L^2_c\to L^2_c, \qquad (Tu)(y)=e^{i\frac b 2y}u(y)
\end{align*}
preserves  $H^1_c$ and satisfies
\begin{align}\label{Equivalence forms}
\mathfrak{a}_b(u,v)&=\mathfrak{\tilde a}_b(Tu,Tv),\quad \forall u,v\in H^1_C,\qquad 
L_b=T^{-1}\circ A_b \circ T.
\end{align}
\end{lem}
{\sc{Proof.}} The proof follows immediately by the equality  $D_yTu=T(D_yu+i\frac b 2 u)$.\qed 

\begin{os} \label{link} It is easier to prove sectoriality and domination for the form $\mathfrak{\tilde a}_b$ and the operator $A_b$ rather than for $\mathfrak {a}_b$ and $L_b$ since $\Rp \mathfrak{\tilde a}_b(u,u)=\|D_y u\|_2^2$ (and this is not true for $\mathfrak{a}_b $).
\end{os}

%
%
%
%
%
%
%

The following results follows easily from  Proposition \ref{kernel Vreal tilde}, since $L_b=T^{-1}\circ A_b \circ T$.

\begin{prop}\label{kernel Vreal}
Let $c+1>0$.   Then the   operator $L_b$  generates an analytic semigroup  of angle $\frac\pi 2$ in $L^2_c$  which satisfies for any $\epsilon>0$

\begin{align*}
	\|e^{z L_b}f\|_{L^2_c}\leq e^{\frac{b^2}{4\epsilon^2 }\Rp z}\|f\|_{L^2_c},\qquad \forall z\in \Sigma_{\frac\pi 2-\arctan \eps }.
\end{align*} 
Moreover
\begin{itemize}
  \item[(i)] The semigroup $\left(e^{t L_b}\right)_{t\geq 0}$ is  $L^\infty$-contractive  and it is dominated by $\left (e^{tB^n}\right )_{t \geq 0}$, that is 
\begin{align*}
 |e^{t L_b}f|\leq e^{tB^n}|f|,\qquad t>0,\quad  f\in L^2_c.
 \end{align*}
   \item[(ii)] $\left(e^{t L_b}\right)_{t\geq 0}$ is a semigroup of integral operators and its heat kernel $p_b$, taken with respect the measure $\rho^cd\rho$, satisfies for some constant $C$ independent of $b$
   \begin{align*}
 |p_{b}(t,y,\rho)|\leq C t^{-\frac{1}{2}}  \rho^{-c}\left (1\wedge\frac{\rho}{\sqrt t} \right)^{c}
\exp\left(-\frac{|y-\rho|^2}{\kappa t}\right)),\qquad \text{for a.e. $y,\rho>0$}.
    \end{align*}     
\item[(iii)] 
For any $s>0$ the operator satisfies the scaling property
\begin{align*}
I_{\frac 1 s}\circ L_b \circ I_s= s^2L_{\frac b s},\qquad I_s u(y):=u(sy). 
\end{align*}
  \end{itemize}  
\end{prop}
{\sc Proof.} All the stated properties follows from 
 Proposition \ref{kernel Vreal tilde} and Lemma \ref{iso} using 
the equalities 
\begin{align*}
	e^{tL_b}f(y)=e^{-i\frac b 2y}\,e^{tA_b}\left(e^{i\frac b 2\,\cdot\, }f\right)(y),\qquad p_b(t,y,\rho)= e^{i\frac b 2(\rho-y)}\,\tilde p_b(t,y,\rho).
\end{align*}
\qed

To  extend the above heat kernel estimates to the half plane $\C_+$ we need the following  elementary lemma, see {\cite[Lemma 5.2]{MNS-PerturbedBessel}} for a straightforward proof.

\begin{lem}\label{Misura palle}
Let $c+1>0$ and for  $y_0,r>0$
\begin{align*}
Q_c (y_0,r):=\int_{y_0}^{y_0+r} y^{c} dy.
\end{align*}
 Then one has 
\begin{align*}
	Q_c(y_0,r)&=r^{c+1}Q_c(\frac{y_0}r,1),\qquad 
Q_c(y_0,r)\simeq r^{c+1}\left(\frac{y_0}{r}\right)^{c}\left(\frac{y_0}{r}\wedge 1\right)^{-c}.
\end{align*}
In particular the function $Q_c$ satisfies, for some  constants $C\geq 1$, the doubling condition
\begin{align*}
	\frac{Q_c(y_0,s)}{Q_c(y_0,r)}\leq C \left(1 \vee \frac{s}r\right)^{1\vee (c+1)},\qquad \forall\, s,  r>0.
\end{align*}
\end{lem}

We also need to rewrite the estimate in Proposition \ref{kernel Vreal}(ii) in an equivalent but more symmetric way.

\begin{lem}\label{os equiv estim}
	  The estimate in  Proposition \ref{kernel Vreal}(ii)  is equivalent (after modifying the constant in the exponential) 
	to 
	\begin{align}\label{equiv estim}
	\nonumber	|p_{b}(t,y,\rho)|&\leq Ct^{-\frac{c+1}2} \left(\frac{y}{\sqrt {t}}\right)^{-\frac c 2} \left(1\wedge \frac{y}{\sqrt {t}}\right)^{\frac c 2}\left(\frac{\rho}{\sqrt {t}}\right)^{-\frac c 2} \left(1\wedge \frac{\rho}{\sqrt {t}}\right)^{\frac c 2}\exp\left(-\frac{|y-\rho|^2}{\kappa t}\right)\\[1.5ex]
		&\simeq \frac{1}{\sqrt{Q_c(y,\sqrt{t}) Q_c(\rho,\sqrt{t})}}\exp\left(-\frac{|y-\rho|^2}{\kappa t}\right).
	\end{align}
\end{lem}
{\sc Proof.} This follows by \cite[Lemma 10.2]{MNS-Caffarelli} with $\gamma_1=\gamma_2=-\frac c 2$.
\qed

\subsection{Bounds for $e^{z L_b}$}
In order to extend the above kernel estimates to complex $z$ we use the standard machinery of \cite[Chapter 6]{Ouhabaz},  relying on estimates like \eqref{equiv estim}  where the terms   $\frac{1}{\sqrt{Q_c(y,\sqrt{t}) Q_c(\rho,\sqrt{t})}} $ are substituted   by powers of $t$, due to the doubling property of Lemma \ref{Misura palle}.

We consider the isometry 
\begin{align} \label{deffi}
	\Phi : L^{2}\left(\R_+,d\mu_{Q_c}\right)\to L^{2}\left(\R_+,d\mu\right), \qquad f\mapsto \frac{f}{\sqrt{Q_{c}(\,\cdot\,,1)}}
\end{align}
where 
\begin{align*}
	d\mu_{Q_c}:=\frac{ d\mu}{Q_{c}(\,\cdot\,,1)}\simeq (y\wedge 1)^c\,dy, 
\end{align*}
by Lemma \ref{Misura palle}. The map $\Phi$ defines a similar operator $\tilde L_b=\Phi^{-1} L_b \Phi $ which acts on $L^{2}\left(\R_+,d\mu_{Q_c}\right)$. We collect in the following proposition the main properties satisfies by $\tilde L_b$ which follow,  by construction,  by   Proposition \ref{kernel Vreal}.

\begin{prop}\label{L tilde main}
$\tilde L_b$ generates an analytic semigroup of angle $\frac \pi 2$ in  $L^{2}\left(\R_+,d\mu_{Q_c}\right)$   which satisfies for any $\epsilon>0$
	\begin{align}\label{tilde L eq 2}
		\|e^{z\tilde L_b}f\|_{L^{2}\left(\R_+,d\mu_{Q_c}\right)}\leq e^{\frac{b^2}{4\epsilon^2 }\Rp z}\|f\|_{L^{2}\left(\R_+,d\mu_{Q_c}\right)},\qquad \forall z\in \Sigma_{\frac\pi 2-\arctan \epsilon }.
	\end{align} 
	Moreover
\begin{itemize}
	  \item[(i)] The semigroup $\left(e^{t\tilde L_b}\right)_{t\geq 0}$ is dominated by $e^{t\tilde L_0}$, that is 
	\begin{align*}
		|e^{t\tilde L_b}f|\leq e^{t\tilde L_0}|f|,\qquad t>0,\quad  f\in L^{2}\left(\R_+,d\mu_{Q_c}\right).
	\end{align*}
\item[(ii)]		
	 $\left(e^{t\tilde L_b}\right)_{t\geq 0}$ is a semigroup of integral operators and its heat kernel $p_{\tilde L_b}$, taken with respect the measure $d\mu_{Q_c}$, satisfies
	\begin{align}\label{nucleo equivalente}
		p_{\tilde L_b}(t,y,\rho)&=\sqrt{Q_{c}(y,1)Q_{c}(\rho,1)}p_{b}(t,y,\rho).
	\end{align}
\end{itemize} 
\end{prop}

\medskip 

The doubling condition of  $Q_c$ guarantees the ultracontractivity of  $\left\{e^{z\tilde L_b}:\ z\in\C_+\right\}$. The following lemma is the main reason for introducing the new operator $\tilde L_b$. Its kernel, in fact, is bounded by a constant depending only on $t$.

\begin{lem}\label{ultracontr tilde}
	$(e^{t\tilde L_b})_{t \geq 0}$ satisfies for some constant $C$ independent of $b$
	\begin{align*}
		& (i)\  \|e^{t\tilde L_b}\|_{\mathcal L\left(L^{1}\left(d\mu_{Q_c}\right),L^{\infty}\right)}\leq C \left(1+ \frac 1 t\right)^{\frac{1\vee (c+1)}2}, \quad (ii)\  
 \|e^{t\tilde L_b}\|_{\mathcal L\left(L^{1}\left(d\mu_{Q_c}\right),L^{2}\left(d\mu_{Q_c}\right)\right)}\leq C \left(1+ \frac 1 t\right)^{\frac{1\vee (c+1)}4} \\[1.5ex]
		 & (iii)\  \|e^{t\tilde L_b}\|_{\mathcal L\left(L^{2}\left(d\mu_{Q_c}\right),L^{\infty}\left(d\mu_{Q_c}\right)\right)}\leq C \left(1+ \frac 1 t\right)^{\frac{1\vee (c+1)}4}.
	\end{align*}
\end{lem}
{\sc Proof.} To prove (i) we observe that  from \eqref{nucleo equivalente},  Lemma \ref{os equiv estim} and Lemma \ref{Misura palle}, one has 
\begin{align}\label{tilde L eq 3}
	\nonumber 	|p_{\tilde L_b}(t,y,\rho)|&=\sqrt{Q_{c}(y,1)Q_{c}(\rho,1)}|p_{b}(t,y,\rho)|
	\leq C \Big[\frac{Q_{c}(y,1)Q_{c}(\rho,1)}{Q_c(y,\sqrt{t}) Q_c(\rho,\sqrt{t})}\Big]^{\frac 1 2}\exp\left(-\frac{|y-\rho|^2}{\kappa t}\right)\\[1.5ex]
	&\leq C \left(1\vee \frac{1}{\sqrt t}\right)^{1\vee (c+1)}\exp\left(-\frac{|y-\rho|^2}{\kappa t}\right)
	\simeq C \left(1+ \frac 1 t\right)^{\frac{1\vee (c+1)}2}\exp\left(-\frac{|y-\rho|^2}{\kappa t}\right).
\end{align}
In particular $e^{t\tilde L_b}$ is ultracontractive and satisfies (i).

To prove (ii) we have for $f\in L^{1}\left(\R_+,d\mu_{Q_c}\right)$, using Proposition \ref{L tilde main}(i) and then (i) of this lemma with $b=0$ (note that $\tilde L_0$ is, by construction,  self-adjoint), 
$$
\|e^{t \tilde L_b}f\|^2_2\leq \|e^{t \tilde L_0}|f| \|^2_2=\langle e^{2t \tilde L_0} |f|, |f| \rangle \leq C \left(1+ \frac 1 t\right)^{\frac{1\vee (c+1)}2}\|f\|_1^2.
$$

Finally,  (iii) follows by duality.\\
\qed 

We can now prove heat kernel estimates for $\tilde L_b$ for complex times. We  need the following  lemma of   Phragm\'en-Lindel\"of type, see \cite[Lemma 6.18]{Ouhabaz}.
\begin{prop}
	\label{Prg-Lind Ouhabaz} Let $\psi\in(0,\frac\pi 2]$ and let $F$ be  an analytic function on $\Sigma_\psi$.
	Assume that, for given numbers $A, \gamma>0$, 
	\begin{alignat}{2}
		\nonumber 
		|F(z)|&\le A(\Rp z)^{-\nu},   &     \qquad\forall& \, z \in \Sigma_\psi ,\\
		\label{Ouh1}
		|F(t)|&\le A t^{-\nu}e^{-\frac{\gamma}{t}},  &   \qquad\forall& \, t>0.
	\end{alignat}	
	Then for any $0<\psi' <\psi$ one has 
	\begin{equation}  \label{Ouh2}
		|F(z)| \le A2^{\nu}(\Rp z)^{-\nu}\exp \left(- \frac{\gamma}{2|z|}\sin(\psi-\psi')\right),
		\quad \forall \, z \in \Sigma_{\psi'} .
	\end{equation}
\end{prop}

\begin{prop}\label{kernel tilde}
	The semigroup $\left\{e^{z\tilde L_b}:\ z\in\C_+\right\}$ consists of integral operators
	\[
	e^{z\tilde L_b}f(y)=
	\int_{0}^\infty p_{\tilde L_b}(z,y,\rho)f(\rho)\,d\mu_{Q_c},\quad f\in L^{2}\left(\R_+,d\mu_{Q_c}\right),\quad y>0.
	\]
	Furthermore for  every $\epsilon>0$, $0<\delta<1$ there exist $C,k>0$ independent of $b$  such that, for every  $z\in\Sigma_{\frac\pi 2-\arctan \epsilon }$ and almost every $y,\rho>0$, one has 
	\begin{align*}
		|p_{\tilde L_b}(z,y,\rho)|\leq C e^{\frac{b^2}{4\epsilon^2\delta}\Rp z} \left(1+ \frac 1 {|z|}\right)^{\frac{1\vee (c+1)}2}\exp\left(-
		\frac{|y-\rho|^2}{k |z|} \right).
	\end{align*}
\end{prop}
{\sc Proof.}  Let us fix  $\epsilon>0$ and  $0<\delta<1$.  Let us  observe that  
\begin{align*}
	t+is\in  \Sigma_{\frac\pi 2-\arctan (\delta \epsilon)}\qquad \Longrightarrow\qquad \delta t+is\in  \Sigma_{\frac\pi 2-\arctan (\delta^2\epsilon) }.
\end{align*}
Then using the semigroup law and  Lemma \ref{ultracontr tilde}   we have for any $z=t+is\in  \Sigma_{\frac\pi 2-\arctan (\delta\epsilon) }$ and for some positive constant $C=C(\epsilon,\delta)$ which may vary in each occurrence
\begin{align}\label{tilde L eq 4}
	\nonumber &\|e^{z\tilde L_b}\|_{\mathcal L\left(L^{1}\left(\R_+,d\mu_{Q_c}\right), L^\infty\left(\R_+\right)\right)}\\[1ex]
	\nonumber &\leq \|e^{\frac{1-\delta }{2}t\tilde L_b}\|_{\mathcal L\left(L^{2}\left(d\mu_{Q_c}\right), L^\infty\right)}\,\|e^{(\delta t+is)\tilde L_b}\|_{\mathcal L\left(L^{2}\left(d\mu_{Q_c}\right), L^2\left(d\mu_{Q_c}\right)\right)}\,\|e^{\frac{1-\delta }{2}t\tilde L_b}\|_{\mathcal L\left(L^{1}\left(d\mu_{Q_c}\right), L^2\left(d\mu_{Q_c}\right)\right)}\\[1ex]
	&\leq C \left(1+ \frac 1{\frac{1-\delta}  2t}\right)^{\frac{1\vee (c+1)}2}e^{\frac{b^2}{4(\delta^2\epsilon)^2 }\delta \Rp z}\simeq C \left(1+ \frac 1{\Rp z}\right)^{\frac{1\vee (c+1)}2}e^{\frac{b^2}{4\epsilon^2\delta^3}\Rp z}
\end{align}

The Dunford-Pettis Theorem then yields the existence of a kernel  $p_{\tilde L_b}(z,y,\rho)$ which satisfies the first claim of the Proposition  and such that 
\begin{align*}
	|p_{\tilde L_b}(z,y,\rho)|\leq C_\epsilon \left(1+ \frac 1{\Rp z}\right)^{\frac{1\vee (c+1)}2}e^{\frac{b^2}{4\epsilon^2\delta^3}\Rp z},\qquad \forall z\in  \Sigma_{\frac\pi 2-\arctan (\delta\epsilon) },\, y,\rho>0.
\end{align*}

Let us now consider the  analytic
function $\Gamma(f_1,f_2):\Sigma_{\frac\pi 2-\arctan (\delta \epsilon)}\rightarrow \mathbb{C}$   defined by

\begin{equation*}
	\Gamma(f_1,f_2)(z)=\langle e^{zL_b}f_{1},f_{2}\rangle_{L^2(d\mu_{Q_c})} ,
\end{equation*}%
where $f_1 \in L^2(F,d\mu_{Q_c})\cap L^1(F,d\mu_{Q_c})$, $f_2\in L^2(E,d\mu_{Q_c})\cap L^1(E,d\mu_{Q_c})$ and $E$, $F$ are two disjoints compact subsets of $\R_+$. Let $r=d(E,F)$ their distance. Then from  \eqref{tilde L eq 3}, \eqref{tilde L eq 4} we have
\begin{alignat*}{2}
	|\Gamma(f_1,f_2)(z)|&\leq C \left(1+ \frac 1{\Rp z}\right)^{\frac{1\vee (c+1)}2}e^{\frac{b^2}{4\epsilon^2\delta^3}\Rp z} \|f_1\|_{L^1(d\mu_{Q_c})}\|f_2\|_{L^1(d\mu_{Q_c})},&\qquad &z\in \Sigma_{\frac\pi 2-\arctan(\delta\epsilon) };\\[1ex]
	|\Gamma(f_1,f_2)(t)|&\le C \left(1+ \frac 1 t\right)^{\frac{1\vee (c+1)}2}\exp\left(-\frac{r^2}{\kappa t}\right)
	\|f_1\|_{L^1(d\mu_{Q_c})}\|f_2\|_{L^1(d\mu_{Q_c})}, & \qquad &t\in\R^+,
\end{alignat*}
(note that $|z|\simeq \Rp z$, $\Rp{\frac 1 z}\simeq \frac 1 {\Rp z}$  for $z\in \Sigma_{\frac\pi 2-\arctan (\delta\epsilon) }$).

Therefore  the function $H(z)=\Gamma(f_1,f_2)(z)e^{-\frac{b^2 z}{4\epsilon^2\delta^3 } }\left(1+\frac 1 z\right)^{-\frac{1\vee (c+1)}2}$ satisfies \eqref{Ouh1} with 
\begin{equation*}
	\gamma =r^{2}/k,\quad A=\|f_1\|_{L^1(d\mu_{Q_c})}\|f_2\|_{L^1(d\mu_{Q_c})},\quad \nu=0.
\end{equation*}%
Proposition \ref{Prg-Lind Ouhabaz} then implies that for some positive constant $C=C(\epsilon,\delta)$, $k=k(\epsilon,\delta)$,  one has  for any $z\in \Sigma_{\frac\pi 2-\arctan \epsilon }$
\begin{align*}
	|\Gamma(f_1,f_2)(z)|&\le C e^{\frac{b^2}{4\epsilon^2\delta^3}\Rp z} \left(1+ \frac 1 {|z|}\right)^{\frac{1\vee (c+1)}2}\exp\left(-\frac{r^2}{\kappa  |z|}\right)
	\|f_1\|_{L^1(d\mu_{Q_c})}\|f_2\|_{L^1(d\mu_{Q_c})}.
\end{align*}

Let us fix $y,\rho>0$  such that  $|y-\rho|>2s>0$ and let us set $r=|y-\rho|-2s$. Then    
\begin{align*}
	|p_{\tilde L_b}(z,y,\rho)|&\leq \sup \left\{	|p_{\tilde L_b}(z,y',\rho')|:\quad y'\in B(y,s),\;\rho'\in
	B(\rho,s) \right\}\\[1ex]
	&= \sup\left\{|\Gamma(f_1,f_2)(z)|\colon \,
	\|f_1\|_{L^1(B(y,s), d\mu_{Q_c})} = \|f_2\|_{L^1(B(\rho,s), d\mu_{Q_c})}=1 \right\} \\[1ex]
	&\leq   C e^{\frac{b^2}{4\epsilon^2\delta^3}\Rp z} \left(1+ \frac 1 {|z|}\right)^{\frac{1\vee (c+1)}2}\exp\left(-\frac{r^2}{\kappa  |z|}\right).
\end{align*}
Recalling that $r=|y-\rho|-2s$ and letting  $s\to 0$ we obtain 
\begin{align*}
	|p_{\tilde L_b}(z,y,\rho)|\leq C e^{\frac{b^2}{4\epsilon^2\delta^3}\Rp z} \left(1+ \frac 1 {|z|}\right)^{\frac{1\vee (c+1)}2}\exp\left(-
	\frac{|y-\rho|^2}{k |z|} \right)
\end{align*}
for any $z\in \Sigma_{\frac\pi 2-\arctan \epsilon }$, which is equivalent to the statement, given the arbitrariness of $\delta$.
\qed

\begin{os}
We remark that in  \cite{Coulhon-Sikora}, the authors work in an abstract metric measure space $(M,d, \mu)$ and assume that the heat kernel $p$  is continuous  with respect to the space variables. In such a case, in fact, 
\begin{align*}
	\underset{x\in U_1, y\in U_2}{ sup}|p(z,x,y)|
	=\mbox{sup}\{\int_{M}e^{-zL}f_1\overline{f_2}\,d\mu,\quad \|f_1\|_{L^1(U_1,d\mu)}=\|f_2\|_{L^1(U_2,d\mu)}=1\}.
\end{align*}
In our setting the continuity assumption  can be avoided since the proofs of \cite[Theorem 4.1, Corollary 4.4]{Coulhon-Sikora}  hold only assuming that  for a.e. $x,y\in M$
\begin{align*}
	p(z,x,y)=\lim_{s\to 0}\int_{M}e^{-zL}f_1\overline{f_2}\,d\mu
	=\lim_{s\to 0} \frac{1}{\mu(B(x,s))\mu(B(y,s))}\int_{B(x,s)\times B(y,s)}p(z,\bar x ,\bar y),d\mu(\bar x)d\mu(\bar y),
\end{align*}
where  $f_1=\frac{\chi_{B(x,s)}}{\mu(B(x,s))}$, $f_2=\frac{\chi_{B(y,s)}}{\mu(B(y,s))}$. 

This holds,  outside a set of zero measure, when the measure $\mu$ is doubling as in our case, by the Lebesgue differentiation theorem.
\end{os}

Finally we  prove  estimates for the heat kernel of $L_b$ for complex times.

\begin{teo}\label{kernel V}
	Let $c+1>0$, $b\in\R$. The semigroup $\left\{e^{zL_b}:\ z\in\C_+\right\}$ consists of integral operators
	\[
	e^{zL_b}f(y)=
	\int_{0}^\infty p_{b}(z,y,\rho)f(\rho)\,\rho^c d\rho,\quad f\in L^2_c,\quad y>0.
	\]
	Furthermore for  every $\epsilon>0$, $0<\delta<1$ there exist $C,k>0$ independent of $b$ such that, for every  $z\in\Sigma_{\frac\pi 2-\arctan \epsilon }$ and almost every $y,\rho>0$, one has 
	\begin{align*}
		|p_b(z,y,\rho)|
		\leq 
		C\,e^{\frac{b^2}{4\epsilon^2\delta}\Rp z}\, |z|^{-\frac{1}{2}}  \rho^{-c}\left (\frac{\rho}{|z|^{\frac{1}{2}}}\wedge 1 \right)^{c}
		\exp\left(-\frac{|y-\rho|^2}{\kappa |z|}\right).
	\end{align*}
\end{teo}
{\sc Proof.} 
The existence of the kernel follows directly using the isometry $\Phi$, see \eqref{deffi}, and  Proposition \ref{kernel tilde}; in particular \eqref{nucleo equivalente}  extends to complex time. Moreover   for some positive constant $C=C(\epsilon,\delta)$   one has 
\begin{align*}
	|p_b(z,y,\rho)|&= \frac{1}{\sqrt{Q_c(y,1)Q_c(\rho,1)}}|p_{\tilde L_b}(z,y,\rho)|\\[2ex]
	&\leq C\, \sqrt{\frac{Q_c(y,\sqrt{|z|})Q_c(\rho,\sqrt{|z|})}{Q_c(y,1)Q_c(\rho,1)}}  \left(1+ \frac 1 {|z|}\right)^{\frac{1\vee (c+1)}2} \\[1ex]
	&\hspace{15ex}\times \frac{1}{\sqrt{Q_c(y,\sqrt{|z|})Q_c(\rho,\sqrt{|z|})}}\,e^{\frac{b^2}{4\epsilon^2\delta}\Rp z}\exp\left(-
	\frac{|y-\rho^2|}{k  |z|} \right).
\end{align*}
Lemma \ref{Misura palle}  then implies
\begin{align}\label{kernel V eq1}
\nonumber 	|p_b(z,y,\rho)|
	&\leq C \left(1+|z|\right)^{\frac{1\vee (c+1)}2}\left(1+ \frac 1 {|z|}\right)^{\frac{1\vee (c+1)}2}\\[1ex]&\hspace{10ex}\times \frac{1}{\sqrt{Q_c(y,\sqrt{|z|})Q_c(\rho,\sqrt{|z|})}}\,e^{\frac{b^2}{4\epsilon^2\delta}\Rp z}\exp\left(-
	\frac{|y-\rho^2|}{k  |z|} \right).
\end{align}
To get rid of the extra term $\left(1+|z|\right)^{\frac{1\vee (c+1)}2}\left(1+ \frac 1 {|z|}\right)^{\frac{1\vee (c+1)}2}$ we use the scaling property (iii) of Proposition \ref{kernel Vreal} and the fact that in the above estimate the constants involved do not depend on $b$. We write for any $z\in\Sigma_{\frac\pi 2-\arctan \epsilon }$, $z=\omega |z|$ and   observe that the scaling equalities of Proposition \ref{kernel Vreal}(iii)
\begin{align*}
	|z|L_b=I_{\frac{1}{\sqrt{|z|}}}\circ  L_{b\sqrt{|z|}} \circ  I_{\sqrt{|z|}},\qquad e^{zL_b}= I_{\frac{1}{\sqrt{|z|}}}\circ e^{\omega L_{b|z|}} \circ I_{\sqrt{|z|}}
\end{align*}
imply
\begin{align*}
 p_b(z,y,\rho)=&|z|^{-\frac{c+1}2}p_{b\sqrt{|z|}}\left(\omega, \frac{y}{\sqrt{|z|}},\frac{\rho }{\sqrt{|z|}} \right).
\end{align*}
Applying \eqref{kernel V eq1} with $b$ replaced by $b\sqrt{|z|}$ to the latter equality  and using Lemma \ref{Misura palle} again we get 
\begin{align*}
	|p_b(z,y,\rho)|&\leq |z|^{-\frac{c+1}2}\left|p_{b\sqrt{|z|}}\left(\omega, \frac{y}{\sqrt{|z|}},\frac{\rho }{\sqrt{|z|}} \right)\right|\\[1ex]
	&\leq C |z|^{-\frac{c+1}2} \frac{1}{\sqrt{Q_c\left(\frac y{\sqrt{|z|}},1\right)Q_c\left(\frac \rho {\sqrt{|z|}},1\right)}}\,e^{\frac{b^2|z|}{4\epsilon^2\delta}\Rp \frac{z}{|z|}}\exp\left(-
	\frac{|y-\rho^2|}{k  |z|} \right)\\[1ex]
	&=C \frac{1}{\sqrt{Q_c(y,\sqrt{|z|})Q_c(\rho,\sqrt{|z|})}}\,e^{\frac{b^2}{4\epsilon^2\delta}\Rp z}\exp\left(-
	\frac{|y-\rho^2|}{k  |z|} \right).
\end{align*}
This concludes the proof, by Lemma \ref{os equiv estim}.
\qed

\subsection{Generation properties and domain characterization}

 First we prove that the semigroup $e^{z L_b}$ extrapolates to the spaces $L^p_m$.

\begin{prop}\label{Gen DN}
	If $c>-1$, $1<p \leq \infty$ and $0< \frac{m+1}{p} <c+1$, then  $(e^{zL_{b}})$ is an  analytic semigroup of angle $\frac\pi 2$ in $L^p_m$ and $\hat{C}([0, \infty[$. Furthermore 
	for every $\epsilon>0$, $0<\delta<1$ there exists $C>0$, independent of $b$,  such that
	\begin{align}\label{dom semigroup}
		\left|e^{zL_b}f\right|\leq Ce^{\frac{b^2}{4\epsilon^2\delta}\Rp z}S^{0,-c}(|z|)|f|,\quad f\in L^p_{m},\quad |\arg z|<\frac\pi 2-\arctan \epsilon 
	\end{align}
	where
	\begin{align*}
		S^{0,-c}(t)f(y)=t^{-\frac {1} 2}\int_{0}^{+\infty}  \left (\frac{\rho}{\sqrt t}\wedge 1 \right)^{c}
		\exp\left(-\frac{|y-\rho|^2}{\kappa t}\right)f(\rho) \,d\rho,
	\end{align*}
	for  a suitable $\kappa>0$.
\end{prop}
{\sc Proof. } All properties for  $p=2$, $m=c$ are contained in Theorem  \ref{kernel V}.  The boundedness of $e^{zL_{b}}$ in $L^p_m$ then  follows from \eqref{dom semigroup} and \cite[Proposition 12.2]{MNS-Caffarelli},  \cite[Proposition 6.2]{MNS-Caf-Schauder} and \eqref{dom semigroup} extends to $L^p_m$.

 The semigroup law is inherited  from the one of $L^2_c$ via a density argument and we have only to prove the  strong continuity at $0$.  Assume first that $p<\infty$ and let $f, g \in C_c^\infty (\R_+)$. Then as $z \to 0$, $z \in \Sigma_{\frac\pi 2-\arctan \epsilon}$, 
$$
\int_0^\infty (e^{zL_{b}}f)\, g\, y^m dy=\int_0^\infty (e^{zL_{b}}f) \,g\, y^{m-c}y^c  dy \to \int_0^\infty fgy^{m-c}y^cdy  =\int_0^\infty fgy^{m}dy,
$$
by the strong continuity of $e^{zL_b}$ in $L^2_c$. Let us observe now that, using  \eqref{dom semigroup} and  \cite[Proposition 12.2]{MNS-Caffarelli},   the family $\left\{e^{z\left(L_{b}-\frac{b^2}{4\epsilon^2\delta }\right)}, z\in \Sigma_{\frac\pi 2-\arctan \epsilon}\right\}$  is uniform bounded on $\mathcal B(L^p_m)$.  Up to rescaling, a   density argument then proves that the previous limit holds for every $f \in L^p_m$, $g \in L^{p'}_m$. The semigroup is then weakly continuous, hence strongly continuous.

When $p=\infty$ the proof  of the strong continuity in $\hat{C}([0,\infty[)$ follows from the domain characterization for $p<\infty$ of Theorem \ref{domainLb} below, as in  \cite[Proposition 2.3]{MNS-Caf-Schauder}.
\qed

To characterize the domain of $L_b$ we need the following interpolation inequality. 

\begin{lem} \label{interpolazione}
If $c>-1$ and $0< \frac{m+1}{p} <c+1$ there exists $C>0$ such that  for every $u \in W^{2,p}_{m, \mathcal N}$
$$
\|D_y u\|_{L^p_m} \leq \eps \|Bu\|_{L^p_m} + \frac C \eps \|u\|_{L^p_m}.
$$
A similar estimate holds for $p=\infty$ if $u \in \hat{C}_{\mathcal N}^2$.
\end{lem}
{\sc Proof. } For $p<\infty$ we use the pointwise estimate
$$
|D_y e^{tB^n} f| \leq \frac{C}{\sqrt t} S_0^{-c}(t)|f|
$$
which follows from Proposition \ref{Estimates Bessel kernels} and the norm estimate $\|S_0^{-c}(t)\|_{L^p_m} \leq C$, see \cite[Proposition 12.2]{MNS-Caffarelli} or Proposition \ref{R-bound S^a,b}. If $u \in  W^{2,p}_{m, \mathcal N}$ and $f=\lambda u-B^n u, \lambda >0$, then $u=\int_0^\infty e^{-\lambda t} e^{t B^n} f\, dt$ and then
\begin{align*}
\|D_y u\|_{m,p} &\leq \int_0^\infty \|e^{-\lambda t}D_y e^{tB^n} f\|_{m,p}\, dt \leq C\int_0^\infty \frac{e^{-\lambda t}}{\sqrt t}\|S_0^{-c}| f|\|_{m,p}\, dt \\
&\leq \frac{C}{\sqrt \lambda} \|f\|_{m,p} \leq C\left(\sqrt \lambda \|u\|_{m,p} +\frac{1}{\sqrt \lambda} \|Bu\|_{m,p}\right).
\end{align*}
When $p=\infty$ the proof is similar, see \cite[Corollary 4.7]{MNS-Caf-Lebesgue}.
\qed
\begin{teo} \label{domainLb}
If $c>-1$, $1<p <\infty$ and $0< \frac{m+1}{p} <c+1$, the generator of $(e^{z L_b})$ is the operator $L_b$ with domain $W^{2,p}_{m,\mathcal N}$. 

When $p=\infty$ the generator is $L_b$ with domain $\hat{C}^2_{\mathcal N}$.
\end{teo}
{\sc Proof. } By the lemma above, the operator $D_y$ is a small perturbation of $B^n$ (with domain $W^{2,p}_{m,\mathcal N}$ or $\hat{C}^2([0, \infty[)$) and  therefore, by \cite[Chapter III, Theorem 2.10]{engel-nagel}, $L_b=B^n+ibD_y$, with the same domain as $B^n$,  generates an analytic semigroup. We have to show that this semigroup coincides what that constructed before. 

We consider first the case $p<\infty$.  
Let $(L_{m,p},D_{m,p})$ be the generator of $(e^{tL_b})$ in $L^p_m$ and consider the set $$\mathcal {D}= \left\{u \in C_c^\infty ([0,\infty)): u \ {\rm constant \ in \ a \ neighborhood\  of }\  0  \right \}$$ which is dense in  $W^{2,p}_{m,\mathcal N}$, by Theorem \ref{Neumnan alpha m 2}. 

By using the definition of $L_b$ through the form  $\mathfrak{a}_b$, see \eqref{defLb}, it is easy to see that $\mathcal D \subset D_{c,2}$ and that $L_b=L_{c,2}$ on $\mathcal D$. Since $\mathcal D$ is dense in $W^{2,2}_{c, \mathcal N}$,  $L_b$ is closed on $W^{2,2}_{c,\mathcal N}$ and $L_{c,2}$ is closed on $D_{c,2}$, it follows that $W^{2,2}_{c,\mathcal N} \subset D_{c,2}$ and then  $W^{2,2}_{c, \mathcal N} = D_{c,2}$,  $L_b=L_{c,2}$ since both operators are generators. This completes the proof in the special case $p=2, m=c$.

Take now $u \in \mathcal D$ and let $f=\lambda u-L_b u \in L^p_m \cap L^2_c$ for large $\lambda$. Let $v \in D_{m,p}$ solve $\lambda v-L_{m,p}v=f$. Since the semigroups are consistent, $v$ coincides with the $L^2_c$ solution which, by the previous step, is $u$. This gives $\mathcal D \subset D_{m,p}$ and that $L_b=L_{m,p}$ on $\mathcal D$ and, as before, one concludes the proof for $p<\infty$.

When $p=\infty$, we change $\mathcal D$  with $$\mathcal 
 D \subset \tilde {\mathcal D}=\{u \in C^2([0, \infty]), D_yu(0)=0, u\ {\rm with\ compact\ support}\} $$ which is a dense in $\hat {C}^2_{\mathcal N}$. We choose $p$ such that $\frac 1p <c+1$ and argue as above, using $L^p$ (with respect to the Lebesgue measure) instead of $L^2_c$.
\qed

Formula \eqref{Equivalence forms} and the previous results allow to  characterize the domain of $A_b=B-i\frac{bc}{2y}$. Note that the Neumann boundary condition for $L_b$, that is  $y^c  v'(y) \to 0$ as  $y \to 0$  translates into a (complex) Robin condition $y^c(u'(y)-i\frac b2 u(y)) \to 0$ for $A_b$. We formulate this result only for $p<\infty$.
\begin{prop}\label{Gen tilde Lb}
	If $c>-1$ and $0< \frac{m+1}{p} <c+1$, then 
	\begin{align*}
		D(A_{b})=\left\{u\in W^{2,p}_m\,:\,\frac 1 y\left(D_yu-i\frac{b}2u\right)\in L^p_m \right\}
	\end{align*}
\end{prop}
{\sc{Proof.}} Let us observe that the isometry $T$ preserves $W^{2,p}_m$ and that  for any $v\in W^{2,p}_m$, setting $u(y)=Tv(y)=e^{i\frac b 2 y}v(y)$,  one has
\begin{align*}
	D_y u=T\left(D_yv+i\frac b 2v\right),\qquad  D_{yy} u=T\left(D_{yy}v+ibD_yv-\frac{b^2}4v \right).
\end{align*} 
 The proof then  follows using the equalities
\begin{align*}
A_b&=T\circ L_b \circ T^{-1},\qquad 	D(A_{b})=T^{-1}\left(W^{2,p}_{m, \mathcal N}\right).
\end{align*}
\qed

\subsection{Bounds for $D_y e^{z L_b}$}

We need the following regularity result which follows from the holomorphy of $e^{z L_b}$ and the characterization of the domain of $L_b$.

\begin{lem} \label{regpb}
For every fixed $\rho>0$ the kernel $p_b(z, y, \rho)$ is holomorphic with respect to $z \in \C_+$ and  $p(z, \cdot , \rho) \in \hat{C}^2_{\mathcal N}$, if $\Rp z >0$.  Moreover all derivatives are jointly continuous in $\C_+ \times [0, \infty[$.
\end{lem}
{\sc Proof.} Fixing $p$ such that $\frac 1p <c+1$ we work in $L^p=L^p_0$. If $s>0$, by Theorem \ref{kernel V}, $p_b(s, \cdot, \rho) \in L^p$ and then $e^{z L_b} p_b$ belongs to the domain of $L_b$  in $L^p$, since the semigroup is analytic. Since $e^{z L_b }p_b(s,y,\rho)=p_b(z+s, y, \rho)$, by the semigroup law, we have that $p(z+s, \cdot, \rho) \in W^{2,p} \subset \hat C([0, \infty[)$. Repeating the argument in this last space we obtain  by Theorem \ref{kernel V} with $p=\infty$ that $p(z+2s, \cdot, \rho) \in \hat{C}^2_{\mathcal N}$.

The analyticity with respect to $ z \in \C_+$ and the joint continuity of the derivatives follow again by the identity $e^{z L_b }p_b(s,y,\rho)=p_b(z+s, y, \rho)$, using the analyticity of the semigroup in $\hat{C}([0, \infty[$, since the domain is $\hat{C}^2_{\mathcal N}$.
\qed

	The Cauchy formula for the derivatives of holomorphic functions allows to estimate $D_zp_b$ and  $L_b p_b$.
	\begin{prop}\label{Time derivative estimates}
			Let $c+1>0$, $b\in\R$. Then  for  every $\epsilon>0$, $0<\delta<1$ there exist $C,k>0$ independent of $b$, such that, for every  $z\in\Sigma_{\frac\pi 2-\arctan \epsilon }$ and almost every $y,\rho>0$,
		\begin{align*} 
			\left|L_b p_b(z,y,\rho)\right|+	\left|D_zp_b (z,y,\rho)\right|
			\leq C  e^{\frac{b^2}{4\epsilon^2\delta}\Rp z} |z|^{-\frac{3}{2}}  \rho^{-c}\left (\frac{\rho}{|z|^{\frac{1}{2}}}\wedge 1 \right)^{c}
			\exp\left(-\frac{|y-\rho|^2}{\kappa |z|}\right).
		\end{align*}
	\end{prop}
{\sc{Proof.}}
		Since the kernel $p_b$  satisfies the equation $D_zp_b=L_b p_b$,   it is sufficient to deal only with $D_z p_b$.   Let us fix $\epsilon>0$, $0<\delta<1$. Setting $r:=\tan\left(\frac{\arctan\epsilon-\arctan(\delta\epsilon)}2\right)<1$, let us observe  that
		\begin{align*}
			B\left(z_0,  r|z_0|\right)\subset \Sigma_{\frac\pi 2-\arctan (\delta\epsilon) }, \qquad \forall z_0\in\Sigma_{\frac\pi 2-\arctan \epsilon }.
		\end{align*}
			Using the Cauchy formula for the derivatives of holomorphic functions in the ball  $B\left(z_0,  r|z_0|\right)$, we get 
		\begin{align*} 
			\left |D_z p_b(z_0,y,\rho)\right|
			\leq \frac{1}{r|z_0|}\max_{|z-z_0|=r|z_0|}|p_b(z,y,\rho)|,\qquad y,\rho>0.
		\end{align*}
		Applying the estimate of  Theorem \ref{kernel V} in the sector $\Sigma_{\frac\pi 2-\arctan (\delta\epsilon) }$, we obtain for suitable $C', \kappa'$
		\begin{align*} 
		\left |D_z p_b(z_0,y,\rho)\right|\leq C'\frac 1{ |z_0|^{\frac 32}}   
			e^{\frac{b^2}{4\epsilon^2\delta^3 }\Rp z_0}   \rho^{-c}\left (\frac{\rho}{|z_0|^{\frac{1}{2}}}\wedge 1 \right)^{c}
			\exp\left(-\frac{|y-\rho|^2}{\kappa' |z_0|}\right)
		\end{align*}
	which, by the arbitrariness of $\delta$,   is equivalent to the statement.
\qed

Before proving the estimates for  the derivative of the kernel $p_b$, let us observe that when $b\neq 0$,  using  the scaling equalities of Proposition \ref{kernel Vreal}(iii), we have   for any $z\in\Sigma_{\frac\pi 2}$
	\begin{align*}
		\frac{1}{b^2}L_b=I_{|b|}\circ  L_{\frac{b}{|b|}} \circ  I_{\frac{1}{|b|}},\qquad e^{zL_b}= I_{|b|}\circ e^{b^2z L_{\frac{b}{|b|}}} \circ  I_{\frac{1}{|b|}}
	\end{align*}
	and then  for  $y, \rho >0$
	\begin{align}\label{Scaling kernel}
	\nonumber	p_b(z,y,\rho)=&|b|^{c+1}\,p_{\frac{b}{|b|}}\left(b^2z,|b|y, |b|\rho \right),\\[2.5ex]
		D_yp_b(z,y,\rho)=&|b|^{c+2}\,D_yp_{\frac{b}{|b|}}\left(b^2z,|b|y, |b|\rho \right).
	\end{align}
	The last equalities allow to assume $|b|=1$ in what follows.

We start by proving  some interpolative estimates with respect to the sup-norm $\|\cdot\|_\infty$.

\begin{lem}\label{lem stime Lb inf}
	Let   $c+1>0$ and $|b|=1$.  Then there exists $C>0$ such that for $\lambda>0$ and  $u\in \hat{C}^2_{\mathcal N}$ one has
	\begin{align*}
		\lambda\|u\|_\infty+\lambda^\frac{1}{2}\|D_y u\|_\infty+\|D_{yy} u\|_\infty\leq C\left( \|\lambda u-L_b u\|_\infty+\|u\|_{\infty}\right).
	\end{align*}
\end{lem} 
%
{\sc Proof.} Since $(L_b,\hat{C}^2_{\mathcal N}) $ generates a bounded semigroup, we have for $u\in \hat{C}^2_{\mathcal N}$, $\lambda >0$
\begin{align*}
	\lambda\|u\|_\infty\leq C\left( \|\lambda u-L_b u\|_\infty\right), \quad \|D_{yy}u\|_\infty \leq C(\|u\|_\infty+\|L_bu\|_\infty).
\end{align*}
Then  
\begin{align*}
	\|D_{yy}u\|_\infty &\leq C\left(\|\lambda u-L_bu\|_\infty+ (\lambda+1) \|u\|_\infty \right)\leq 2C\left( \|\lambda u-L_bu\|_\infty+\|u\|_{\infty}\right).
\end{align*}

The estimate of the gradient term follows by the inequality $\|D_y u\|_\infty^2 \leq C\|u\|_\infty \|D_{yy}u\|_\infty$.
\qed

Now we localize the gradient estimates above. For   $y>0$, $r>0$ we set $B^+(y,r):=B(y,r)\cap\R_+$. 

	\begin{prop}\label{Proposition inner estimates} Let $c+1>0$, $|b|=1$. Then there exists a constant $C>0$  such that  for every $u\in \hat{C}_{\mathcal N}$ and  $\lambda>0$ 
	\begin{align*} 
		\lambda^{\frac 12} \|D_yu\|&_{L^\infty(B^+(y,r))}\leq C\left(\|\lambda u-L_bu\|_{L^\infty(B^+(y,2r))}+\left(\frac{1}{ r^2}+1\right)\|u\|_{L^\infty(B^+(y,2r))}\right).
	\end{align*} 
\end{prop}
{\sc Proof.}  Set $r_n=r\sum_{k=1}^n2^{-k}$. Then $r_1=r/2$, $r_\infty=r$, $r_{n+1}-r_n=r2^{-(n+1)}$.

Let  $B^+_n=B^+(y,r_n)$ $B_r^+= B^+(y,r)$ and so on, and choose  cut-off functions $\eta_n\in C_c^{\infty }(\R)$ such that  $\eta_n(y)=\eta_n(-y)$,
$0\le \eta_n\le 1$, 
$\eta_n=1$ in $B_n^+$, $({\rm supp \ } \eta_n) \cap \R_+ \subset B_{n+1}^+$,
$|D_y \eta_n| \leq \frac{C}{r}2^n$,  $|D_{yy} \eta_n| \leq \frac{C}{r^2}4^n$ for some constant $C>0$ independent of $n$. Then also $|y^{-1}D_y \eta_n| \leq \frac{C}{r^2}4^n$, since $D_y \eta_n(0)=0$.  

If  $u\in \hat{C}_{\mathcal N}$ then  $\eta_n u\in \hat{C}_{\mathcal N}$   and we have
$$L_b(\eta_n u)=\eta_n L_bu +2 D_y \eta_n D_y u +u\left( D_{yy}\eta_n+c\frac{D_y\eta_n}{y}+ibD_y\eta_n\right).$$
Applying Lemma \ref{lem stime Lb inf} to $\eta_n u$ and using the  inequality $s\leq 1+s^2$  we get 
\begin{align*}
	\lambda \|\eta_nu\|_\infty+&\sqrt \lambda \|D_y (\eta_nu)\|_\infty+\|D_{yy} (\eta_nu)\|_\infty\leq C\left( \|(\lambda-L_b )(\eta_nu)\|_\infty+ \|\eta_nu\|_\infty\right)\\[1.5ex]
	&\leq C\left( \|\lambda u-L_b u\|_{L^\infty(B_r^+)}+\frac{2^n}r\|D_yu\|_{L^\infty(B_{r_{n+1}}^+)}\right.\\[1.5ex] &\Bigg.\hspace{10ex}+\|u\|_{L^\infty(B_r^+)}\left\|D_{yy}\eta_n+c\frac{D_y\eta_n}{y}+ibD_y\eta_n\right\|_{\infty}+\|u\|_{L^\infty(B_r^+)}\Bigg)\\[1.5ex]
	&\leq C\left( \|\lambda u-L_b u\|_{L^\infty(B_r^+)}+\frac{2^n}r\|D_y(\eta_{n+1}u)\|_{\infty}+\left(\frac{4^n}{r^2}+1\right)\|u\|_{L^\infty(B_r^+)}\right).
\end{align*}
Applying the interpolative inequalities for the gradient we get
\begin{align*}
	\lambda \|\eta_nu\|_\infty+&\sqrt\lambda\|D_y (\eta_nu)\|_\infty+\|D_{yy} (\eta_nu)\|_\infty\leq C\left( \|\lambda u-L_b u\|_{L^\infty(B_r^+)}+\epsilon \frac{2^n}r\|D_{yy}(\eta_{n+1}u)\|_{\infty}\right.\\[1.5ex]
	&\left.\hspace{28ex}+\frac{2^n}{\epsilon r}\|\eta_{n+1}u\|_{\infty}+\left(\frac{4^n}{r^2}+1\right)\|u\|_{L^\infty(B_r^+)}\right)\\[1.5ex]
	&\leq C\left( \|\lambda u-L_b u\|_{L^\infty(B_r^+)}+\epsilon \frac{2^n}r\|D_{yy}(\eta_{n+1}u)\|_{\infty}+\left(\frac{2^n}{\epsilon r}+\frac{4^n}{r^2}+1\right)\|u\|_{L^\infty(B_r^+)}\right)
\end{align*}
Setting
$\xi:=C2^n\eps r^{-1}$, we get
\begin{align*}
	\lambda \|\eta_nu\|_\infty+&\sqrt\lambda\|D_y (\eta_nu)\|_\infty+\|D_{yy} (\eta_nu)\|_\infty\\[1.5ex]
	&\leq C\left( \|\lambda u-L_b u\|_{L^\infty(B_r^+)}+\left(\frac{4^n}{\xi r^2}+\frac{4^n}{r^2}+1\right)\|u\|_{L^\infty(B_r^+)}\right)+\xi \|D_{yy}(\eta_{n+1}u)\|_{\infty}
\end{align*}
It follows that 
\begin{align*}
	&\xi^n \left(\sqrt\lambda \|D_y u \|_{L^\infty(B_{\frac r 2}^+)}+\|D_{yy} (\eta_nu)\|_\infty\right)\\[1.5ex]
	&\leq C\left( \xi^n\|\lambda-L_b u\|_{L^\infty(B_r^+)}+\xi^n\left(\frac{4^n}{\xi r^2}+\frac{4^n}{r^2}+1\right)\|u\|_{L^\infty(B_r^+)}\right)+\xi^{n+1} \|D_{yy}(\eta_{n+1}u)\|_{\infty}
\end{align*}
By choosing $\eps=\eps_n$  so that $\xi=\frac 18$ and summing up previous inequality over $n\in\N$ we get
\begin{align*}
	& \sqrt\lambda \|D_y u \|_{L^\infty(B_{\frac r 2}^+)}+\sum_{n=1}^\infty \xi^n \|D_{yy} (\eta_nu)\|_\infty\\[1.5ex]
	&\leq C\left(\|\lambda-L_b u\|_{L^\infty(B_r^+)}+\left(\frac{1}{ r^2}+1\right)\|u\|_{L^\infty(B_r^+)}\right)+\sum_{n=1}^\infty \xi^{n+1} \|D_{yy}(\eta_{n+1}u)\|_{\infty}
\end{align*}
Cancelling equal terms on both sides it follows that
\begin{align*}
	\sqrt\lambda \|D_y u \|_{L^\infty(B_{\frac r 2}^+)}\leq C\left(\|\lambda u-L_b u\|_{L^\infty(B_r^+)}+\left(\frac{1}{ r^2}+1\right)\|u\|_{\infty,r}\right).
\end{align*}
\qed

To prove estimates for  the derivative of the kernel $p_b$ we use also the following basic estimate.
\begin{lem}\label{lemma loc gauss}
	Let $y_0, \rho >0$. Then 
	\begin{align*}
\sup_{y\in B(y_0,\sqrt t) }\exp\left(-\frac{|y-\rho|^2}t\right)\leq e^{16}\exp\left(-\frac 9{16}\frac{|y_0-\rho|^2}t\right)
	\end{align*}
\end{lem}
{\sc{Proof.}}  If $|y_0-\rho|\leq 4\sqrt t$ then for every $y >0$
\begin{align*}
\exp\left(-\frac{|y-\rho|^2}t\right)\leq  e^{16}\exp\left(-\frac{|y_0-\rho|^2}t\right).
\end{align*}
If $|y_0-\rho|> 4\sqrt t$ and $y\in B(y_0,\sqrt t)$. then $|y-\rho|\geq |y_0-\rho|-|y-y_0| \geq |y_0-\rho|-\frac 1 4 |y_0-\rho|=\frac 3 4 |y_0-\rho|$ and 
\begin{align*}
	\exp\left(-\frac{|y-\rho|^2}t\right)\leq  \exp\left(-\frac 9{16}\frac{|y_0-\rho|^2}t\right).
\end{align*}
\qed 
\begin{teo}\label{space derivative estimates}
	Let $c+1>0$, $b\in\R$.  Then  for  every $\epsilon>0$, $0<\delta<1$ there exist $C,k>0$ independent of $b$, such that, for every  $z\in\Sigma_{\frac\pi 2-\arctan \epsilon }$ and almost every $y,\rho>0$,
	\begin{align*} 
			\left|D_yp_b (z,y,\rho)\right|
		\leq C e^{\frac{b^2}{4\epsilon^2\delta}\Rp z} \frac1{|z|}  \rho^{-c}\left (\frac{\rho}{|z|^{\frac{1}{2}}}\wedge 1 \right)^{c}
		\exp\left(-\frac{|y-\rho|^2}{\kappa |z|}\right).
	\end{align*}
\end{teo}
{\sc{Proof.}} If $b=0$ this follows from  Proposition \ref{Estimates Bessel kernels}. Let now $b\neq 0$; using the scaling property \eqref{Scaling kernel}, we may  assume that $|b|=1$.  Then applying Proposition \ref{Proposition inner estimates}  to the function  $u=p_b(z,\,\cdot\,,\rho)$ in  $B^+(y_0,r)$  with $r=\sqrt{|z|}$  we get for any $\lambda>0$ 
\begin{align*} 
		\lambda \left|D_yp_b (z,y_0,\rho)\right|\leq C\left(\|L_bu\|_{L^\infty(B^+(y,2r))}+\left(\frac{1}{r^2}+1+\lambda^2\right)\|u\|_{L^\infty(B^+(y,2r))}\right).
\end{align*} 
 Using Theorem \ref{kernel V} and  Proposition \ref{Time derivative estimates} with $\delta'$ such that $0<\delta<\delta'<1$,  we get for suitable  $C, \kappa >0$ 
  \begin{align*} 
 	\|u\|_{L^\infty(B^+(y_0,2r))}&\leq  C e^{\frac{1}{4\epsilon^2\delta'}\Rp z} |z|^{-\frac{1}{2}}  \rho^{-c}\left (\frac{\rho}{|z|^{\frac{1}{2}}}\wedge 1 \right)^{c}
 	\sup_{y\in B^+(y_0,2r) }\exp\left(-\frac{|y-\rho|^2}{\kappa |z|}\right),\\[2ex]
 	\|L_bu\|_{L^\infty(B^+(y_0,2r))}&\leq  C e^{\frac{1}{4\epsilon^2\delta'}\Rp z} |z|^{-\frac{3}{2}}  \rho^{-c}\left (\frac{\rho}{|z|^{\frac{1}{2}}}\wedge 1 \right)^{c}
 	\sup_{y\in B^+(y_0,2r) }\exp\left(-\frac{|y-\rho|^2}{\kappa |z|}\right).
 \end{align*}
Lemma \ref{lemma loc gauss} then implies (for suitable $C', \kappa'>0$)
 \begin{align*} 
 	\|u\|_{L^\infty(B^+(y_0,2r))}\leq A(z,y_0,\rho),\qquad 	\|L_bu\|_{L^\infty(B^+(y_0,2r))}\leq  \frac{1}{|z|}A(z,y_0,\rho),\\[2ex]
 	A(z,y_0,\rho):= C' e^{\frac{1}{4\epsilon^2 \delta'}\Rp z} |z|^{-\frac{1}{2}}  \rho^{-c}\left (\frac{\rho}{|z|^{\frac{1}{2}}}\wedge 1 \right)^{c}
 	\exp\left(-\frac{|y_0-\rho|^2}{\kappa' |z|}\right).
 \end{align*}
 The previous inequalities then imply 
\begin{align*} 
	\left|D_yp_b (z,y_0,\rho)\right|&\leq \left(\frac{1}{\lambda 
	}\left(\frac{2}{ |z|}+1\right)+\lambda\right)A(z,y_0,\rho).
\end{align*}
Minimizing the last inequality over $\lambda>0$ we get for a suitable $C>0$
\begin{align*} 
	\left|D_yp_b (z,y_0,\rho)\right|
	\leq 2\frac{\sqrt{2+|z|}}{\sqrt{|z|}}\;A(z,y_0,\rho)\leq C\frac{e^{\frac{1}{4\epsilon^2 }\left(\frac{1}{\delta}-\frac{1}{\delta'}\right)\Rp z}}{\sqrt{|z|}}\;A(z,y_0,\rho).
\end{align*}
which is the statement for $|b|=1$. 
\qed 

\begin{cor}\label{Dy T(t)_b}
Let  $1<p \leq \infty$, $0< \frac{m+1}{p} <c+1$, $f\in L^p_{m}$ or $f\in\hat{C}([0, \infty[)$.  For  every $\epsilon>0$, $0<\delta<1$ there exist $C>0$ independent of $b$, such that, for every  $z\in\Sigma_{\frac\pi 2-\arctan \epsilon }$ ,   $e^{zL_{b}}f$ is differentiable in $]0, \infty[$ and satisfies
	\begin{align}\label{der semigroup}
		& D_ye^{zL_b}f=\int_0^{+\infty}D_yp_b(z,\cdot,\rho)f(\rho)\,\rho^{c} d\rho\\ \label{der semigroup1}
		& |D_ye^{zL_b}f|\leq \frac{C}{ |z|^{\frac{1}{2}}}e^{\frac{b^2}{4\epsilon^2\delta}\Rp z} S^{0,-c}(|z|)|f|.
	\end{align}
	\end{cor}
{\sc Proof.} 
	Let $y_0, r>0$ such that $0\not\in[y_0-r, y_0+r]$. By Theorem \ref{space derivative estimates},  for almost every $y\in ]y_0-r, y_0+r[$, $\rho\in]0,\infty[$,
\begin{align}  \label{stimagradiente}
			\left|D_yp_b (z,y,\rho)\right|\rho^c
		\leq C e^{\frac{b^2}{4\epsilon^2\delta}\Rp z} \frac1{ |z|}  \left (\frac{\rho}{|z|^{\frac{1}{2}}}\wedge 1 \right)^{c}
		\exp\left(-\frac{\rho^2}{\kappa |z|}\right)
	\end{align}
	for  suitable  $C$ and $ \kappa$ depending also  on $r$ and $y_0$. Then \eqref{der semigroup} follows by differentiating under the integral sign since 
the right hand side of \eqref{stimagradiente} belongs to $L^{p'}_m$. Finally, \eqref{der semigroup1} is consequence of Theorem \ref{space derivative estimates}.
	\qed

\section{Multipliers} \label{mult}

In this section we investigate the boundedness of some multipliers related to  the   degenerate operator
	\begin{align*}
		\mathcal L =\Delta_{x} +2\sum_{i=1}^Na_{i}D_{iy}+D_{yy}+\frac{c}{y}D_y ,\qquad a \in\R^N,\ |a|<1,
	\end{align*}
	Assuming that
	$$\Delta_x u+2a\cdot \nabla_xD_yu+ B_yu=f$$
	and taking the Fourier transform (denoted by $\mathcal F$ or $\hat \cdot$)  with respect to $x$ (with covariable $\xi$) we obtain 
	$$ -|\xi|^2\,\hat u(\xi,y)+i2a\cdot\xi D_y\hat u(\xi,y)+  B _y \hat u(\xi,y)=\hat f(\xi,y).$$ 
We  introduce the quadratic form 
\begin{align}\label{est Qa}
	Q_{a}(\xi)=|\xi|^2-|a\cdot\xi|^2,\qquad(1-|a|^2)|\xi|^2\leq Q_{a}(\xi)\leq |\xi|^2,\qquad  \xi\in\R^N
\end{align}
and  consider the operator $L_b$ of Section \ref{section L_b} with $b=2a\cdot\xi$. The latter computation  shows that formally 
\begin{align*}
	(\lambda-\mathcal L)^{-1}={\cal F}^{-1}\left (\lambda+Q_a(\xi)-L_{2a\cdot\xi}  \right)^{-1} {\cal F}.
\end{align*} 
In order to prove that $\mathcal L$ generates an analytic semigroup and to prove regularity for the associated parabolic problem, we   investigate  the boundedness of the  operator-valued multiplier
 \begin{align*}
 	\xi \in \R^N \to R_{\lambda}(\xi)&= (\lambda+Q_a(\xi)-L_{2a\cdot\xi})^{-1}.
 \end{align*}
 To characterize the  domain of $\mathcal L$ we also consider the multipliers   $|\xi|^2R_\lambda$, $\xi D_yR_{\lambda}$ which are associated respectively with the operators $	\Delta_x(\lambda-\mathcal L)^{-1}$, $ D_{xy}	(\lambda-\mathcal L)^{-1}$. In the next results  we  prove that the  above  multipliers  satisfy the hypotheses of Theorem \ref{marcinkiewicz}. 

 \smallskip

 The following lemma is a reformulation of the heat kernel bounds of the previous section, adapted to the multipliers above.
For any  $|a|\leq  \delta<1$ we  set 
 \begin{align*}
 	\theta_\delta=\arctan\frac{|a|}{\sqrt{\delta^2-|a|^2}}\in\left]\arctan\frac{|a|}{\sqrt{1-|a|^2}},\frac\pi 2\right]
 \end{align*}
where $\theta_{|a|}:=\frac\pi 2$.

\begin{lem}\label{R-bound Semigroup}
	Let  $1<p<\infty$ be such that $0<\frac{m+1}p<c+1$.  For every $|a|\leq \delta<\delta'< 1$, there exists $C>0$ such that for $ f\in L^p_{m}$, $  z \in \Sigma_{\frac\pi 2-\theta_\delta}$
		\begin{align*}
			\left|e^{z(L_{2a\cdot\xi}-Q_a(\xi))}f\right|\leq C  e^{-\left(Q_{\delta' \frac{a}{|a|}}(\xi)\right)\Rp z}S^{0,-c}(|z|)|f|
			\leq C e^{-(1-\delta'^2)|\xi|^2\Rp z}S^{0,-c}(|z|)|f|
		\end{align*}
		and
		\begin{align*}
			\left|D_ye^{z(L_{2a\cdot\xi}-Q_a(\xi))}f\right|&\leq C  e^{-(1-\delta'^2)|\xi|^2\Rp z}\frac{S^{0,-c}(|z|)}{\sqrt{|z|}}|f|.
		\end{align*}
\end{lem}  
{\sc Proof. } 
We use Proposition  \ref{Gen DN} with  $\eps=\frac{|a|}{\sqrt{\delta^2-|a|^2}}$, $0<\gamma=\frac{\delta^2-|a|^2}{\delta'^2-|a|^2}<1$ and \eqref{est Qa} which yield
\begin{align*}
	\left|e^{z(L_{2a\cdot\xi}-Q_a(\xi))}f\right|&\leq C e^{-\left(Q_a(\xi)-\frac{|a\cdot\xi|^2}{\epsilon^2\gamma}\right)\Rp z}S^{0,-c}(|z|)|f|\\[1ex]
	&=C e^{-\left(Q_{\delta' \frac{a}{|a|}}(\xi)\right)\Rp z}S^{0,-c}(|z|)|f|\leq C e^{-(1-\delta'^2)|\xi|^2\Rp z}S^{0,-c}(|z|)|f|
\end{align*}
for any $f\in L^p_{m}$, $|\arg z|<\frac\pi 2-\theta_\delta$. The estimate for $D_ye^{z(L_{2a\cdot\xi}-Q_a(\xi))}$ follows similarly using Theorem \ref{space derivative estimates}.
\qed

The following formulas follow since the resolvent is the Laplace transform of the semigroup.  We state them to have precise bounds and to show that a similar representation holds for the gradient of the resolvent.

\begin{lem}\label{Explicit Resolv}
	Let   $|a|\leq  \delta<\delta'<1$.	 The following properties hold for any $f\in L^p_m$ and  $\lambda\in\Sigma_{\pi-\theta_\delta}$
		\begin{align*}
			\left(\lambda+Q_a(\xi)-L_{2a\cdot\xi}\right)^{-1}=e^{-i\theta}\int_0^\infty e^{-e^{-i\theta}\lambda  t}e^{ e^{-i\theta}t(L_{2a\cdot\xi}-Q_a(\xi))}f\,dt
		\end{align*}
		and 
		\begin{align*}
			D_y\left(\lambda+Q_a(\xi)-L_{2a\cdot\xi}\right)^{-1}&=e^{-i\theta}\int_0^\infty e^{-e^{-i\theta}\lambda  t}D_ye^{ e^{-i\theta}t(L_{2a\cdot\xi}-Q_a(\xi))}f\,dt
		\end{align*}
Here $\theta=\frac{|\mbox{arg}\lambda|}{\mbox{arg}\lambda}\left(\frac\pi 2-\theta_{\delta  '}\right)$.
\end{lem}
{\sc{Proof.}} 
	If  $\lambda \in \Sigma_{\pi-\theta_\delta}$,  $\mu:=e^{-i\theta}\lambda \in\Sigma_{\frac\pi 2-\theta_\delta+\theta_{\delta'}}$. By  Lemma \ref{R-bound Semigroup}
	\begin{align*}
		\left|e^{t\left(e^{-i\theta}(L_{2a\cdot\xi}-Q_z(\xi))\right)}f\right|\leq C e^{-(1-\delta')\cos\theta t}S^{0,-c}(t)|f|.
	\end{align*}
This implies, by   standard results on analytic semigroups,   
	\begin{align*}
		\left(\lambda+Q_a(\xi)-L_{2a\cdot\xi}\right)^{-1}&=e^{-i\theta}\Big(e^{-i\theta}\lambda-e^{-i\theta}(L_{2a\cdot\xi}-Q_a(\xi))\Big)^{-1}=\\[1.5ex]
		&=e^{-i\theta}\int_0^\infty e^{-e^{-i\theta}\lambda  t}e^{ e^{-i\theta}t(L_{2a\cdot\xi}-Q_a(\xi))}f\,dt
	\end{align*}
The remaining equality follows similarly by differentiating under the integral (note that both integrals converge).\qed
\begin{teo}\label{R-boundednessMult}
	Let  $1<p <\infty$ be such that $0<\frac{m+1}p<c+1$.  For every $|a|\leq  \delta<1$   there exist $C, k>0$ such that for $f\in L^p_m$ and $\lambda\in\Sigma_{\pi-\theta_\delta}$
		\begin{align}\label{dom resolvent eq1}
			\left|\left(\lambda+Q_a(\xi)-L_{2a\cdot\xi}\right)^{-1}f\right|\leq C\, \Gamma\Big(k(|\lambda|+|\xi|^2)\Big)|f|,
		\end{align}
	\begin{align}\label{dom D resolvent eq2}
		\left|D_y\left(\lambda+Q_a(\xi)-L_{2a\cdot\xi}\right)^{-1}f\right|\leq C\, \Psi\Big(k(|\lambda|+|\xi|^2)\Big)|f|.
	\end{align}
\end{teo}  
{\sc Proof. } 
	For $\lambda \in \Sigma_{\pi-\theta_\delta}$ let us choose $\delta_2<\delta'_2$ such that 
\begin{align*}
|a|\leq  \delta<\delta_2<\delta'_2<1,\qquad  	\theta_\delta-\theta_{\delta_2}\leq \frac\pi 2
\end{align*} and let us set  $\theta=\frac{|\mbox{arg}\lambda|}{\mbox{arg}\lambda}\left(\frac\pi 2-\theta_{\delta  _2}\right)$ so that $\mu:=e^{-i\theta}\lambda \in\Sigma_{\frac\pi 2-\theta_\delta+\theta_{\delta_2}}$. Then  using  Lemma \ref{R-bound Semigroup} and Lemma \ref{Explicit Resolv}  we get for some $C$ depending on $\delta_2$
	\begin{align*}
&	\left|\left(\lambda+Q_a(\xi)-L_{2a\cdot\xi}\right)^{-1}f\right|=  \left|\left(\mu-e^{-i\theta}(L_{2a\cdot\xi}-Q_a(\xi))\right)^{-1}f\right|= \left|\int_0^\infty e^{-\mu t}e^{e^{-i\theta} t(L_{2a\cdot\xi}-Q_a(\xi))}f\,dt\right|\\[1ex]
	&\leq C  \int_0^\infty  e^{-\Rp\mu t}e^{-Q_{\delta'_2\frac{a}{|a|}}(\xi)\cos\theta t}S^{0,-c}(t)|f|\,dt
\leq C\int_0^\infty  e^{-|\lambda|\sin(\theta_{\delta}-\theta_{\delta_2}) t}e^{-Q_{\delta'_2\frac{a}{|a|}}(\xi)\sin\theta_{\delta_2} t}S^{0,-c}(t)|f|\,dt.
	\end{align*}
Then using \eqref{est Qa} we get for some positive constant  $k$ depending on $\delta$
\begin{align*}
	\left|\left(\lambda+Q_a(\xi)-L_{2a\cdot\xi}\right)^{-1}f\right|& \leq C\int_0^\infty  e^{-k\left(|\lambda|+|\xi|^2\right)}S^{0,-c}(t)|f|\,dt=C \Gamma\left(k(|\lambda|+|\xi|^2)\right)|f|,
\end{align*}
which is \eqref{dom resolvent eq1}. The proof of \eqref{dom D resolvent eq2} is similar: using Lemma \ref{R-bound Semigroup} and Lemma \ref{Explicit Resolv} and proceeding as before we have 
	\begin{align*}
	&\left|D_y\left(\lambda+Q_a(\xi)-L_{2a\cdot\xi}\right)^{-1}f\right|=  \left|\left(\mu-e^{-i\theta}D_y(L_{2a\cdot\xi}-Q_a(\xi))\right)^{-1}f\right|\\[1.5ex]
	&\leq C\int_0^\infty  e^{-k_\delta\left(|\lambda|+|\xi|^2\right)}S^{0,-c}(t)|f|\,dt=C \Psi\left(k(|\lambda|+|\xi|^2)\right)|f|.
\end{align*}
\qed

The following result is a consequence of Lemma  \ref{R-bound Semigroup} and Theorem \ref{R-boundednessMult}.
\begin{cor}
	Let  $1<p<\infty$ be such that $0<\frac{m+1}p<c+1$. For every $|a|\leq  \delta<1$  
		the families of operators  
		\begin{align*}
			&\left\{ e^{z(L_{2a\cdot\xi}-Q_a(\xi))}  :\; z\in \Sigma_{\frac \pi 2-\theta_\delta},\xi\in\R^N\setminus\{0\}\right\} \\[1ex]
&			\left\{\lambda\left(\lambda-L_{2a\cdot\xi}+Q_a(\xi)\right)^{-1}, \quad \sqrt{\lambda}D_y\left(\lambda+Q_a(\xi)-L_{2a\cdot\xi}\right)^{-1}:\; \lambda \in \Sigma_{\pi-\theta_\delta},\ \xi\in\R^N\setminus\{0\}\right\}
		\end{align*}
		are $\mathcal R$-bounded in $L^p_m$. 
\end{cor}
{\sc{Proof.}}  The $\mathcal R$-boundedness of $e^{z(L_{2a\cdot\xi}-Q_a(\xi))}$ follows from Lemma \ref{R-bound Semigroup} and  by domination using Corollary \ref{domination}, since  the family $(S^{\alpha,-c}(t))_{t\geq 0}$ is $\mathcal R$-bounded by  Proposition \ref{R-bound S^a,b}(i). 
The $\mathcal R$-boundedness of the families involving the resolvent  follows from domination again,   using Proposition \ref{R-bound S^a,b}(ii) and (iii) since by \eqref{dom resolvent eq1} one has
\begin{align*}
	\left|\lambda\left(\lambda+Q_a(\xi)-L_{2a\cdot\xi}\right)^{-1}f\right|\leq C |\lambda|\,\Gamma\left(k(|\lambda|+|\xi|^2)\right)|f|\leq C |\lambda|\,\Gamma\left(k|\lambda|\right)|f|
\end{align*}
and similarly for $\sqrt{\lambda}D_y\left(\lambda+Q_a(\xi)-L_{2a\cdot\xi}\right)^{-1}$.
\qed

From now on we denote by $R_{\lambda}(\xi)$ the operator $(\lambda-L_{2a\cdot\xi}+Q_a(\xi))^{-1}$ whenever it is defined.

To apply the Mikhlin multiplier theorem, we need a formula for the derivatives  of the above functions with respect to $\xi$. In the following lemma ${\mathcal S}_n$ denotes the set of permutations of $n$ elements.

\begin{lem}   \label{Lema Marck1}
	Let $1<p<\infty$  be such that $0<\frac{m+1}p<c+1$, and let   us consider, for any fixed $\lambda\in\Sigma_{\pi-\arctan \frac{|a|}{\sqrt{1-|a|^2}}}$, the  map  
	\begin{align*}
		\xi \in \R^N \to R_{\lambda}(\xi)&= (\lambda-L_{2a\cdot\xi}+Q_a(\xi))^{-1}\in B(L^p_m).
	\end{align*} 
Then $R_\lambda,D_yR_\lambda\in C^\infty\left(\R^N\setminus \{0\}; B(L^p_m)\right)$ and for any multi-index $\alpha=(\alpha_1,\dots,\alpha_N)\in\{0,1\}^N$, $|\alpha|=n$ one has
\begin{align}\label{Explicit Derivative Rlamba}
\nonumber	D_\xi^\alpha R_\lambda(\xi)&=\sum_{\sigma\in {\mathcal S}_n} R_\lambda(\xi) \prod_{j=1}^n\Big( 2ia_{\sigma(j)}D_yR_\lambda(\xi)-2\xi_{\sigma(j)}R_\lambda(\xi)\Big)\\[1ex]
	D_\xi^\alpha D_yR_\lambda(\xi)&=\sum_{\sigma\in {\mathcal S}_n} D_yR_\lambda(\xi) \prod_{j=1}^n\Big( 2ia_{\sigma(j)}D_yR_\lambda(\xi)-2\xi_{\sigma(j)}R_\lambda(\xi)\Big).
\end{align} 
Furthermore for  every $|a|\leq \delta<1$   there exists $C,k >0$, depending only on  $\delta$  such that, setting $\mu=k(|\lambda|+|\xi|^2)$ one has 
  \begin{align}\label{est Derivative Rlamba}
 \nonumber 	\left|D_\xi^\alpha R_\lambda(\xi)f\right|&\leq C \Gamma\left(\mu\right) \Big( \Psi(\mu)+|\xi|\Gamma\left(\mu\right)\Big)^n|f|, \qquad f\in L^p_m,\quad \lambda\in\Sigma_{\pi-\theta_\delta}\\[1ex]
 	\left|D_\xi^\alpha D_yR_\lambda(\xi)f\right|&\leq C \Psi\left(\mu\right) \Big( \Psi(\mu)+|\xi|\Gamma\left(\mu\right)\Big)^n|f|, \qquad f\in L^p_m,\quad \lambda\in\Sigma_{\pi-\theta_\delta}
  \end{align} 
\end{lem}
{\sc{Proof.}}
Let us fix $|a|\leq \delta<1$ and $\lambda\in\Sigma_{\pi-\arctan \frac{|a|}{\sqrt{\delta-|a|^2}}}$. Let us prove the first equality in  \eqref{Explicit Derivative Rlamba} for $n=1$  that is, for $j=1,\dots, n$   
\begin{align}\label{Lema Marck1 eq2}
	\frac{\partial}{\partial \xi_j} (R_{\lambda}(\xi))= R_{\lambda}(\xi) \left(2ia_jD_y-2\xi_j\right) R_{\lambda}(\xi),\qquad \xi\in\R^n\setminus\{0\}.
\end{align}
Indeed let us  write for $|h|\leq 1$
\begin{align}\label{Lema Marck1 eq1}
	\nonumber  R_{\lambda}&(\xi+he_j)-R_{\lambda}(\xi)= \left(\lambda+Q_a(\xi+he_j)-L_{2a\cdot (\xi+he_j)}\right)^{-1}-\left(\lambda+Q_a(\xi)-L_{2a\cdot \xi}\right)^{-1}\\[2ex]
	\nonumber&= R_{\lambda}(\xi)\,\Big[\left(\lambda+Q_a(\xi)-L_{2a\cdot \xi}\right)\left(\lambda+Q_a(\xi+he_j)-L_{2a\cdot (\xi+he_j)}\right)^{-1}-I\Big]\\[2ex]
	\nonumber&= R_{\lambda}(\xi)\,\left(L_{2a\cdot (\xi+he_j)}-L_{2a\cdot \xi}+Q_a(\xi)-Q_a(\xi+he_j)\right)\left(\lambda+Q_a(\xi+he_j)-L_{2a\cdot (\xi+he_j)}\right)^{-1}\\[2ex]
\nonumber	&=R_{\lambda}(\xi)\,\left(2ia\cdot he_j D_y+|\xi|^2-|\xi+he_j|^2\right)R_{\lambda}(\xi+he_j)\\[2ex]
	&=R_{\lambda}(\xi)\,\left(2ia_jhD_y-2\xi_jh-h^2\right)R_{\lambda}(\xi+he_j).
\end{align}
Applying the previously inequality we get
\begin{align*}
	\frac{R_{\lambda}(\xi+he_j)-R_{\lambda}(\xi)}h&-R_{\lambda}(\xi) \left(2ia_jD_y-2\xi_j\right) R_{\lambda}(\xi)=\\[2ex]
	&=R_{\lambda}(\xi)\,\left(2ia_jD_y-2\xi_j\right)\Big(R_{\lambda}(\xi+he_j)-R_{\lambda}(\xi)\Big)-hR_{\lambda}(\xi)R_{\lambda}(\xi+he_j)\\[2ex]
	&:=A+F.
\end{align*}
From now on we write  $C, k$  to denote some positive constants which depends only on $\delta $.

The term $F$ tends to $0$  as $h$ tends to $0$ in the  norm of $B\left(L^p_m\right)$ since   by Theorem \ref{R-boundednessMult} we have  any $f\in\L^p_m$
	\begin{align*}
	\left|Ff\right|&\leq C\,|h|\, \Gamma\Big(k(|\lambda|+|\xi|^2)\Big)\Gamma\Big(k(|\lambda|+|\xi+he_j|^2)\Big)|f|\leq C\,|h|\,\Gamma\Big(k(|\lambda|)\Big)^2|f|
\end{align*}
where in the last inequality we used the fact that  the  $\Gamma(\lambda)f$ is decreasing in $\lambda>0$ if $f \geq 0$.   This shows that  $\|F\|_{L^p_m}\leq C |h|$. Concerning  $A$, we  apply again \eqref{Lema Marck1 eq1} to write 
\begin{align*}
	A&=R_{\lambda}(\xi)\,\left(2ia_jD_y-2\xi_j\right)R_{\lambda}(\xi)\,\left(2ia_jhD_y-2\xi_jh-h^2\right)R_{\lambda}(\xi+he_j).
\end{align*}
Using again  Theorem \ref{R-boundednessMult} and since  $\Gamma(\lambda)f,\Psi  (\lambda )f$ are decreasing in $\lambda>0$ if $f \geq 0$ we obtain
\begin{align*}
	|Af|&\leq C\, h^2\left|R_\lambda(\xi)D_yR_\lambda(\xi)D_yR_\lambda(\xi+he_j)f\right|+|h||\xi|\left|R^2_\lambda(\xi)D_yR_\lambda(\xi+he_j)f\right|\\[1.5ex]
	&+C\left (h^2(|\xi|+|h|)\left|R_\lambda(\xi)D_yR_\lambda(\xi)R_\lambda(\xi+he_j)f\right|+|h||\xi|(|\xi|+|h|)\left|R^2_\lambda(\xi)R_\lambda(\xi+he_j)f\right|\right )\\[1.5ex]
	\leq &C \,\Big[h^2\,\Gamma(k|\lambda|)\left( \Psi(k|\lambda|)\right)^2+|\xi||h| (\Gamma(k|\lambda|))^2 \Psi(k|\lambda|)\Big.\\[1.5ex]
	&\Big.\hspace{4ex}+|h|^2\left(|\xi|+|h|\right)\Gamma(k|\lambda|)\Psi(k|\lambda|)\Gamma(k|\lambda|)+|h||\xi|\left(|\xi|+|h|\right)(\Gamma(k|\lambda|))^3\Big]|f|.
\end{align*}
Since $\Gamma (k|\lambda|), \Psi(k|\lambda|)$ are bounded in $L^p_m$, by Proposition \ref{R-bound S^a,b}, the last inequality shows that $A$ tends to $0$ in the norm of $B\left(L^p_m\right)$ as $h\to 0$. This proves \eqref{Lema Marck1 eq2}.

The proof of  the second equality in  \eqref{Explicit Derivative Rlamba} for $n=1$  that is, for $j=1,\dots, n$ ,  
\begin{align}\label{Lema Marck1 eq2 D_y}
	\frac{\partial}{\partial \xi_j} (D_yR_{\lambda}(\xi))= D_yR_{\lambda}(\xi) \left(2ia_jD_y-2\xi_j\right) R_{\lambda}(\xi),\qquad \xi\in\R^n\setminus\{0\}
\end{align}
follows similarly and we only sketch the main steps. As in \eqref{Lema Marck1 eq1}  we write
\begin{align*}
	&\frac{D_yR_{\lambda}(\xi+he_j)-D_yR_{\lambda}(\xi)}h-D_yR_{\lambda}(\xi) \left(2ia_jD_y-2\xi_j\right) R_{\lambda}(\xi)=\\[2ex]
	&=D_yR_{\lambda}(\xi)\,\left(2ia_jD_y-2\xi_j\right)\Big(R_{\lambda}(\xi+he_j)-R_{\lambda}(\xi)\Big)
-hD_yR_{\lambda}(\xi)R_{\lambda}(\xi+he_j)
:=A+F.
\end{align*}
Proceeding as before we then  get   for $f\in\L^p_m$ 
	\begin{align*}
	\left|Ff\right|\leq& C\,|h|\,\Psi(k|\lambda|)\Gamma(k|\lambda|)|f|,\\[1ex]
	|Af|\leq& C \,\Big[h^2\,\left( \Psi(k|\lambda|)\right)^3+|h||\xi| \Psi(k|\lambda|)\Gamma(k|\lambda|) \Psi(k|\lambda|)\Big.\\[1.5ex]
	&\Big.\hspace{4ex}+|h|^2\left(|\xi|+|h|\right)\Psi(k|\lambda|)^2\Gamma(k|\lambda|)+|h||\xi|\left(|\xi|+|h|\right)\Psi(k|\lambda|)(\Gamma(k|\lambda|))^2\Big]|f|
\end{align*}
and conclude as before. This proves \eqref{Lema Marck1 eq2 D_y}.

In particular from \eqref{Lema Marck1 eq2} and \eqref{Lema Marck1 eq2 D_y},  using again Theorem \ref{R-boundednessMult}, we get 
\begin{align*}
	\left|\frac{\partial}{\partial \xi_j} (R_{\lambda}(\xi))f\right|&\leq C\left( \Gamma(\mu) \Psi(\mu)+  |\xi|(\Gamma(\mu))^2 \right ),\\[1ex]
	\left|\frac{\partial}{\partial \xi_j} (D_yR_{\lambda}(\xi))f\right|&\leq C\left ( (\Psi(\mu))^2 +  |\xi|\Psi(\mu)\Gamma(\mu) \right ),
\end{align*}
which is \eqref{est Derivative Rlamba} for $n=1$.

Finally,  \eqref{Explicit Derivative Rlamba} for $n>1$  follows by induction. For example if $n=2$ and $l\neq j$ one has 
\begin{align*}
	\frac{\partial^2}{\partial \xi_l\partial \xi_j} (R_{\lambda}(\xi))&=\frac{\partial}{\partial \xi_l}\Big[ R_{\lambda}(\xi) \left(2ia_jDy-2\xi_j\right) R_{\lambda}(\xi)\Big]\\[2ex]
	&=\frac{\partial}{\partial \xi_l}(R_{\lambda}(\xi))\left(2ia_jDy-2\xi_j\right) R_{\lambda}(\xi)+ R_{\lambda}(\xi)\left(2ia_jDy-2\xi_j\right)\frac{\partial}{\partial \xi_l}(R_{\lambda}(\xi))\\[2ex]
	&=R_{\lambda}(\xi) \left(2ia_lDy-2\xi_l\right) R_{\lambda}(\xi)\left(2ia_jDy-2\xi_j\right) R_{\lambda}(\xi)\\[2ex]
	&\hspace{2ex}+ R_{\lambda}(\xi)\left(2ia_jDy-2\xi_j\right) 
	R_{\lambda}(\xi) \left(2ia_lDy-2\xi_l\right) R_{\lambda}(\xi).
\end{align*}

The estimates  \eqref{est Derivative Rlamba} now follow  using again  \eqref{Explicit Derivative Rlamba} and Theorem \ref{R-boundednessMult}.
\qed
Now we can finally prove that the multiplier  $\lambda R_\lambda(\xi)$ associated with the operators  $\lambda(\lambda-\mathcal L)^{-1}$ satisfies the hypothesis of  Theorem \ref{marcinkiewicz}. This is crucial    for proving   that $	\mathcal L $ generates an analytic semigroup  in $L^p_m$.

\begin{teo}\label{Teo Mult Resolv}
		Let $1<p<\infty$ be  such that $0<\frac{m+1}p<c+1$.
	Then for  every $|a|\leq \delta<1$  the family 
	$$\left \{\xi^{\alpha}D^\alpha_\xi(\lambda  R_{\lambda})(\xi): \xi\in \R^{N}\setminus\{0\}, \ \alpha  \in \{0,1\}^N ,\lambda \in \Sigma_{\pi-\theta_\delta}\right \}$$
	is $\mathcal{R}$-bounded in $L^p_m$.
\end{teo}
{\sc{Proof.}} Let $\alpha=(\alpha_1,\dots,\alpha_N)\in\{0,1\}^N$, $|\alpha|=n$ and  $|a|\leq \delta<1$. Using  \eqref{est Derivative Rlamba} we get for some positive constant $C,k>0$ and for any $f\in L^p_m,\quad \lambda\in\Sigma_{\pi-\theta_\delta}$
\begin{align}\label{est Derivative Rlamba 2}
	\left|D_\xi^\alpha R_\lambda(\xi)f\right|&\leq C \Gamma\left(\mu\right) \Big( \Psi(\mu)+\sqrt \mu \Gamma\left(\mu\right)\Big)^n|f|, 
\end{align} 
where $\mu=k\left(|\lambda|+|\xi|^2\right)$. In particular

\begin{align*}
	\left|\xi^{\alpha}\lambda D_\xi^\alpha R_\lambda(\xi)f\right|&\leq C \mu^{\frac {n}2}  \mu\,\Gamma\left(\mu\right) \left( \Psi\left(\mu\right)+\sqrt\mu\, \Gamma\left(\mu\right)\right)^n|f|
= C\mu \Gamma\left(\mu\right) \Big( \sqrt \mu\Psi\left(\mu\right)+\mu\, \Gamma\left(\mu\right)\Big)^n|f|.
\end{align*}
 The  $\mathcal R$-boundedness of $\xi^{\alpha}D^\alpha_\xi(\lambda  R_{\lambda})(\xi)$ then follows by composition and domination using  Proposition \ref{R-bound S^a,b}   and  Corollary \ref{domination}.
 \qed

The next two theorems  show  that the multipliers  $|\xi|^2 R_\lambda$, $\xi D_yR_\lambda$ associated with the operators  $\Delta_x(\lambda-\mathcal L)^{-1}$, $D_{xy}(\lambda-\mathcal L)^{-1}$ satisfy the hypotheses of   Theorem \ref{marcinkiewicz}. This is essential     for characterizing the domain of  $\mathcal L $.
\begin{teo}\label{Mult xi^2 Resolv}
	Let $1<p<\infty$ be  such that $0<\frac{m+1}p<c+1$. 
	Then for  every $|a|\leq \delta<1$  the family 
	$$\left \{\xi^{\alpha}D^\alpha_\xi(|\xi|^2  R_{\lambda})(\xi): \xi\in \R^{N}\setminus\{0\}, \ \alpha  \in \{0,1\}^N ,\lambda \in \Sigma_{\pi-\theta_\delta}\right \}$$
	is $\mathcal{R}$-bounded in $L^p_m$.
\end{teo}
{\sc{Proof.}} Let us observe that for any  $\alpha  \in \{0,1\}^N$, $|\alpha|=n$  there exist $\beta^j  \in \{0,1\}^N$, $|\beta^j|=n-1$, such that
\begin{align*}
D^\alpha_\xi(|\xi|^2  R_{\lambda})(\xi)=\sum_{j:\alpha_j=1}2\xi_jD^{\beta^j}_\xi R_{\lambda}(\xi)+|\xi|^2D^\alpha_\xi R_{\lambda}(\xi). 
\end{align*}
Using \eqref{est Derivative Rlamba} and proceeding as in Theorem \ref{Teo Mult Resolv}, the equality above implies that  for some $C,k>0$ and any $f\in L^p_m$ 
\begin{align*}
	\left|\xi^\alpha D^\alpha_\xi(|\xi|^2  R_{\lambda})(\xi)f\right|\leq& C \mu^{\frac {n+1} 2}\, \Gamma\left(\mu\right) \Big( \Psi\left(\mu\right)+\sqrt\mu\, \Gamma\left(\mu\right)\Big)^{n-1}|f| \\
&+C\mu^{\frac {n} 2+1}\,\Gamma\left(\mu\right) \Big( \Psi\left(\mu\right)+\sqrt\mu\, \Gamma\left(\mu\right)\Big)^{n}|f|
\end{align*} 
with $\mu=k\left(|\lambda|+|\xi|^2\right)$.
The proof now follows as at  the end  of Theorem \ref{Teo Mult Resolv}.\qed

\begin{teo}\label{Mult D_y/y Resolv}
	Let $1<p<\infty$ be  such that $0<\frac{m+1}p<c+1$. 
	Then for  every $|a|\leq \delta<1$  the family 
	$$\left \{\xi^{\alpha}D^\alpha_\xi\left( \xi D_yR_{\lambda}\right)(\xi): \xi\in \R^{N}\setminus\{0\}, \ \alpha  \in \{0,1\}^N ,\lambda \in \Sigma_{\pi-\theta_\delta}\right \}$$
	is $\mathcal{R}$-bounded in $L^p_m$.
\end{teo}
{\sc{Proof.}}  Let us fix $j=1,\dots,N$ and let   $\alpha  \in \{0,1\}^N$, $|\alpha|=n$; then  there exists $\beta  \in \{0,1\}^N$, $|\beta|=n-1$ such that
\begin{align*}
	D^\alpha_\xi(\xi_{j} D_yR_{\lambda})(\xi)=\xi_{j}D^{\alpha}_\xi D_yR_{\lambda}(\xi)+\alpha_{j}D^\beta_\xi D_yR_{\lambda}(\xi). 
\end{align*}
Using \eqref{est Derivative Rlamba} and proceeding as in Theorem \ref{Teo Mult Resolv},  we get   for some  $C,k>0$ and any $f\in L^p_m$ 
\begin{align*}
	\left|\xi^\alpha D^\alpha_\xi(\xi_{j}  D_yR_{\lambda})(\xi)f\right|&\leq C \left(\mu^{\frac {n+1} 2}\, \Psi(\mu) \Big( \Psi\left(\mu\right)+\sqrt\mu\, \Gamma\left(\mu\right)\Big)^{n}|f|	+\mu^{\frac {n} 2}\,\Psi\left(\mu\right) \Big( \Psi\left(\mu\right)+\sqrt\mu\, \Gamma\left(\mu\right)\Big)^{n-1}|f|\right)\\[1ex]
	=& C \left(\sqrt\mu\, \Psi(\mu) \Big( \sqrt\mu\Psi\left(\mu\right)+\mu\, \Gamma\left(\mu\right)\Big)^{n}|f|
+\sqrt\mu\,\Psi\left(\mu\right) \Big( \sqrt\mu\Phi\left(\mu\right)+\mu\, \Gamma\left(\mu\right)\Big)^{n-1}|f|\right)
\end{align*}
where $\mu=k\left(|\lambda|+|\xi|^2\right)$ and the rest is similar.
\qed

\section{Domain and maximal regularity for $\mathcal L=\Delta_{x} +2a\cdot\nabla_xD_y+ B_yu$} \label{Sect DOm max}

In this section we prove generation results, maximal regularity and   domain characterization for  the  degenerate   operator 
 %
\begin{equation}  \label{La}
	\mathcal L :=\Delta_{x} +2\sum_{i=1}^Na_{i}D_{x_iy}+D_{yy}+\frac{c}{y}D_y=\Delta_x u+2a\cdot \nabla_xD_yu+ B_yu ,\qquad a=(a_1, \dots, a_N) \in\R^N,\ |a|<1
\end{equation}
in $L^p_m$. The condition $|a| <1$ is equivalent to the ellipticity of the top order coefficients. More general operators will be treated in the next section, based on this model case. We start with the $L^2$ theory.

\subsection{The operator $\mathcal L$ in $L^2_{c}$} 
We use the Sobolev space $H^{1}_{c}:=\{u \in L^2_{c} : \nabla u \in L^2_{c}\}$ equipped with the inner product
\begin{align*}
	\left\langle u, v\right\rangle_{H^1_{c}}:= \left\langle u, v\right\rangle_{L^2_{c}}+\left\langle \nabla u, \nabla v\right\rangle_{L^2_c}.
\end{align*}
We consider the  form  in $L^2_{c}$ 
\begin{align*}
	\mathfrak{a}(u,v)
	&:=
	\int_{\R^{N+1}_+} \langle \nabla u, \nabla \overline{v}\rangle\,y^{c} dx\,dy+2\int_{\R^{N+1}_+} D_yu\, a\cdot\nabla_x\overline{v}\,y^{c} dx\,dy , \quad
	D(\mathfrak{a})=H^1_{c}
\end{align*}

and its adjoint $\mathfrak{a}^*(u,v)=\overline{\mathfrak{a}(v,u)}$  
\begin{align*}
	\mathfrak{a^*}(u,v)=\overline{\mathfrak{a}(v,u)}
	&:=
	\int_{\R^{N+1}_+} \langle \nabla u, \nabla \overline{v}\rangle\,y^{c} dx\,dy+2\int_{\R^{N+1}_+} a\cdot \nabla_x u\,D_y\overline{v}\,y^{c} dx\,dy.
\end{align*}

\begin{prop} \label{prop-form}
The forms $\mathfrak{a}$, $\mathfrak{a^*}$ are continuous, accretive and   sectorial.
\end{prop}
{\sc Proof.} We consider only the form $\mathfrak{a}$, the adjoint form can be handled similarly. If $u\in H^1_{c}$
\begin{align*}
\Rp \mathfrak{a}(u,u)\geq \|\nabla_x u\|^2_{L^2_c}+\|D_y u\|^2_{L^2_c}-2|a| \|\nabla_x u\|_{L^2_c}\|D_y u\|_{L^2_c}\geq (1-|a|)(\|\nabla_x u\|^2_{L^2_c}+\|D_y u\|^2_{L^2_c}).
\end{align*}
By the ellipticity assumption $|a|<1$, the accretivity follows.
Moreover 
\begin{align*}
|\Ip \mathfrak{a}(u,u)|\leq  2|a|\|\nabla_x u\|_{L^2_c}\|D_y u\|_{L^2_c}\leq |a|(\|\nabla_x u\|^2_{L^2_c}+\|D_y u\|^2_{L^2_c})\leq \frac{|a|}{(1-|a|)}\Rp \mathfrak{a}(u,u).
\end{align*}
This proves the sectoriality and then the continuity of the form. \qed

We define the operators $\mathcal L$ and $\mathcal L^*$ associated respectively to the forms $\mathfrak{a}$ and $\mathfrak{a}^*$  by
\begin{align} \label{BesselN}
	\nonumber D( \mathcal L)&=\{u \in H^1_{c}: \exists  f \in L^2_{c} \ {\rm such\ that}\  \mathfrak{a}(u,v)=\int_{\R^{N+1}_+} f \overline{v}y^{c}\, dz\ {\rm for\ every}\ v\in H^1_{c}\},\\  \mathcal Lu&=-f;
\end{align}
\begin{align} \label{adjoint}
	\nonumber D( \mathcal L^*)&=\{u \in H^1_{c}: \exists  f \in L^2_{c} \ {\rm such\ that}\  \mathfrak{a}^*(u,v)=\int_{\R^{N+1}_+} f \overline{v}y^{c}\, dz\ {\rm for\ every}\ v\in H^1_{c}\},\\  \mathcal L^*u&=-f.
\end{align}
If $u,v$ are smooth function with compact support in the closure of $\R_+^{N+1}$ (so that they do not need to vanish on the boundary), it is easy to see integrating by parts that $$-\mathfrak a (u,v)= \langle \Delta_x u+2a\cdot \nabla_xD_yu+ B_yu, \overline v\rangle_{L^2_c}  $$
if $\lim_{y \to 0} y^c D_yu(x,y)=0$. This means that 
$\mathcal L$ is the operator $\Delta_x +2a\cdot \nabla_xD_y+ B_y$ with Neumann boundary conditions at $y=0$. On the other hand 
$$-\mathfrak a^* (u,v)= \langle \Delta_x u+2a\cdot \nabla_xD_yu+2c\frac {a\cdot \nabla_x u}{y}+ B_yu, \overline v\rangle_{L^2_c}  $$
if $\lim_{y \to 0} y^c \left ( D_yu(x,y)+2a\cdot \nabla_x u(x,y)\right )=0$ and therefore $\mathcal L^*$ is the operator $\Delta_x +2a\cdot \nabla_xD_y+2c\frac {a\cdot \nabla_x u}{y}+ B_y$ with the above oblique condition at $y=0$. 

\begin{prop}\label{generation L2}
 $\mathcal L$ and $\mathcal L^*$  generate contractive analytic semigroups $e^{z \mathcal L}$, $e^{z \mathcal L^*}$, $z\in\Sigma_{\frac{\pi}{2}-\arctan \frac{|a|}{1-|a|}}$,  in $L^2_{c}$. Moreover the semigroups $(e^{t\mathcal L})_{t \geq 0}, (e^{t\mathcal L^*})_{t \geq 0}$ are positive and $L^p_c$-contractive for $1 \leq p \leq \infty$.
\end{prop}
{\sc Proof.} We argue only for $\mathcal L$.
The  generation  result immediately follows from   Proposition \ref{prop-form} and \cite[Theorem 1.52]{Ouhabaz}.
The positivity follows by \cite[Theorem 2.6]{Ouhabaz} after observing that, if $u\in H^1_{c}$, $u$ real, then $u^+\in H^1_{c}$ and
\begin{align*}
	\mathfrak{a}(u^+,u^-)
	&:=
	\int_{\R^{N+1}_+} \langle \nabla u^+, \nabla u^-\rangle\,y^{c} dx\,dy+2\int_{\R^{N+1}_+}D_yu^+ a\cdot \nabla_xu^-\,y^{c} dx\,dy=0.
\end{align*}
Finally, the $L^\infty$-contractivity follows by \cite[Theorem 2.13]{Ouhabaz}  after observing that if $u\in H^1_{c}$,  then $(|u|-1)^+\sign u\in H^1_{c}$ and
\begin{align*}
\Rp \mathfrak{a}(u,(|u|-1)^+\sign u)\geq& \|\nabla_x u\|^2_{L^2_c\cap\{|u|\geq 1\}}+\|D_y u\|^2_{L^2_c\cap\{|u|\geq 1\}}\\[1ex]
&-2|a| |\nabla_x u\|_{L^2_c\cap\{|u|\geq 1\}}\|D_y u\|_{L^2_c\cap\{|u|\geq 1\}}\\[1ex]\geq& (1-|a|)(\|\nabla_x u\|^2_{L^2_c\cap\{|u|\geq 1\}}+\|D_y u\|^2_{L^2_c\cap\{|u|\geq 1\}})\geq 0.
\end{align*}
\qed

The Stein interpolation theorem then shows that the above semigroups are analytic in $L^p_c$ for $1<p<\infty$, see \cite[Proposition 3.12]{Ouhabaz} and a result by Lamberton  yields maximal regularity in the same range, see \cite[Theorem 5.6]{KunstWeis}. Since our results are more general, we do not state these consequences here.

Our aim is to characterize the  domain of $ \mathcal L$ in $L^2_c$ but first we identify a core.

\begin{prop}\label{core Neumnan comp 2}
	If   $c+1>0$   then 
	the set $C_c^\infty (\R^{N})\otimes \mathcal D$ defined in \eqref{defDcore1}
	is a core for $\mathcal L$ in $L^2_{c}$.
\end{prop}
{\sc{Proof.}}  We observe, preliminarily, that  the set $C_c^\infty (\R^{N})\otimes \mathcal D$ is contained in $H^1_{c}$. Moreover, integrating by parts one sees that any $u\in C_c^\infty (\R^{N})\otimes \mathcal D$ satisfies  \eqref{BesselN} with $\mathcal L u=\Delta_x u+2a\cdot\nabla_xD_y u+B_y u \in L^2_{c}$. This yields  $C_c^\infty (\R^{N})\otimes \mathcal D\subseteq D( \mathcal L )$.

Since $I-  \mathcal{L}$ is invertible we have to show that $(I-  \mathcal{L})\left(C_c^\infty (\R^{N})\otimes \mathcal D\right)$ is dense in $L^2_{c}$ or, equivalently,  that $\left((I- \mathcal{L})\left(C_c^\infty (\R^{N})\otimes \mathcal D\right)\right )^{\perp}=\left\{0\right\}$. Let $v\in L^2_{c}$ be such that
\begin{align*}
	\int_{\R^{N+1}_+}\left(I- \mathcal{L}\right)f\, \bar{v}\ dx\ y^{c} dy
	=0, \quad \forall f\in C_c^\infty (\R^{N})\otimes \mathcal D.
\end{align*}
Let us choose $f=a(x)u(y)\in C_c^\infty (\R^{N})\otimes\mathcal{D}$.  Taking the  Fourier transform with respect to $x$   we get $\hat f(\xi,y)=\hat a(\xi) u(y)$ and
\begin{align*}
	\int_{\R^{N+1}_+}\Big[u(y)+ |\xi|^2 u(y)-2ia\cdot\xi D_y u(y)- B_y u(y)\Big]\,\hat a(\xi) \ \bar{\hat v}(\xi,y)\ d\xi\, y^{c}dy=0,
\end{align*}
that is 
\begin{align}\label{Core2 eq 1}
	\int_{\R^{N+1}_+}\Big[u(y)+ Q_a(\xi) u(y)-L_{2a\cdot\xi }\, u(y)\Big]\,\hat a(\xi) \ \bar{\hat v}(\xi,y)\ d\xi\, y^{c}dy=0.
\end{align}
Fix $\xi_0\in\R^N$, $r>0$ and let $w(\xi)=\frac{1}{|B(\xi_0,r)|}\chi_{B(\xi_0,r)}\in L^2(\R^N)$. Let  $(a_n)_n\in C_c^\infty(\R^{N})$ a sequence of test functions such that $a_n\to \check{w}$ in $L^2(\R^N)$; then  $\hat a_n\to w$ in $L^2(\R^N)$ and writing  \eqref{Core2 eq 1} with $\hat a$ replaced by $\hat a_n$ and letting  $n\to\infty$ we obtain
\begin{align*}
	\frac{1}{|B(\xi_0,r)|}\int_{B(\xi_0,r)}d\xi\int_{0}^\infty \Big[u(y)+ Q_a(\xi) u(y)-L_{2a\cdot\xi }\, u(y)\Big]\,\bar{\hat v}(\xi,y)\ y^{c}dy=0.
\end{align*}

Letting $r\to 0$ and using the Lebesgue differentiation theorem, we have for a.e. $\xi_0\in\R^N$
\begin{align*}
	\int_{0}^\infty \left(I+ Q_a(\xi_0)-L_{2a\cdot\xi_0 }\right) u(y)\,\bar{\hat v}(\xi_0,y)\ y^{c}dy=0,\qquad \forall u\in \mathcal{D}.
\end{align*}
 Since, by Theorems  \ref{core gen} and \ref{domainLb},
$\mathcal{D}$ is a core for the operator $L_{2a\cdot\xi_0 }$  in $L^2_{c}(\R_+)$, the last equation  then implies $\hat v(\xi_0,\cdot)=0$  for a.e. $\xi_0\in\R^N$ and the proof is complete.
\qed

\begin{teo}\label{Neumnan comp 2}
	If  $c+1>0$  then 
	\begin{align*}
		D(\mathcal L)&=W^{2,2}_{c,\mathcal{N}}
	\end{align*}
\end{teo}
{\sc Proof. } Observe that
\begin{equation*} 
	C_c^\infty (\R^{N})\otimes \mathcal D\subset W^{2,2}_{c,\mathcal{N}} \cap D(\mathcal L)
\end{equation*}
is a core for $\mathcal L$ by Proposition \ref{core Neumnan comp 2} and is dense in  $W^{2,2}_{c,\mathcal{N}}$ by Theorem \ref{core gen}.

We have to show that the graph norm and that of $W^{2,2}_{c,\mathcal{N}}$ are equivalent on $C_c^\infty (\R^{N})\otimes \mathcal D$. Since the second is obviously stronger, we have to show the converse.

We use Proposition \ref{Trace D_yu in W} and  endow  $W^{2,2}_{c,\mathcal{N}}$ with the equivalent  norm $$\|u\|_W=\|u\|_{L^2_{c}}+\|\Delta_x u\|_{L^2_{c}}+\|\nabla_xD_yu\|_{L^2_{c}}+\| B_yu\|_{L^2_{c}}.$$

Let $u \in  C_c^\infty (\R^{N})\otimes \mathcal D$ and $f= u-\mathcal L u$, so that  $\|u\|_{L^2_{c}} \le \|f\|_{L^2_{c}}$.     By taking the Fourier transform with respect to $x$ (with co-variable $\xi$) we obtain
\begin{align*}(1+Q_a(\xi)-L_{2a\cdot\xi})\hat u(\xi,\cdot)=\hat f(\xi,\cdot)
\end{align*}
and therefore 
\begin{align}\label{Four eq}
\nonumber	|\xi|^2 \hat u(\xi, \cdot)&=|\xi|^2(1+Q_a(\xi)-L_{2a\cdot\xi})^{-1}\hat f(\xi,\cdot)\\[1ex]
	\xi_i D_y\hat u(\xi, \cdot)&=\xi_i D_y(1+Q_a(\xi)-L_{2a\cdot\xi})^{-1}\hat f(\xi,\cdot)
\end{align}

This means that
\begin{align*}
\Delta_x u=-{\cal F}^{-1} |\xi|^2R_1(\xi) {\cal F} f,\qquad  i\frac{\partial}{\partial x_i}D_y u={\cal F}^{-1} \xi D_yR_1(\xi) {\cal F} f.
\end{align*}

The estimates $\| \Delta_x u\|_{L^2_{c}} \le C\|f\|_{L^2_{c}}$, $\|\nabla_xD_yu\|_{L^2_{c}} \le C\|f\|_{L^2_{c}}$  then follow  from the  boundedness of the multipliers $|\xi|^2R_1(\xi)$ and $\xi D_yR_1(\xi)$ in $L^2(\R^N; L^2_{c}(\R_+))=L^2_c$ which are proved in Theorems \ref{Mult xi^2 Resolv} and 
\ref{Mult D_y/y Resolv} and yield  $\| B_y u\|_{L^2_{c}} \le C\|f\|_{L^2_{c}}$,  by difference. 

This gives  the equivalence of the graph norm and of the norm of $W^{2,2}_{c,\mathcal{N}}$ on $C_c^\infty (\R^{N})\otimes \mathcal D$ and concludes the proof.
\qed
\subsection{The operator $\mathcal L$ in $L^p_m$}

In this section we prove  domain characterization and maximal regularity for  $\mathcal L$ 
in $L^p_m$. For clarity reasons we often write $ \mathcal{L}_{m,p}$ to emphasize the underlying space on which the operator acts.

We shall use extensively the set $\mathcal{D}$ defined in \eqref{defDcore}. In particular $\mathcal L$ is well defined on $C_c^\infty (\R^{N})\otimes\mathcal D$ when $m+1 >0$. 

In   Proposition \ref{generation L2} and Theorem \ref{Neumnan comp 2} we proved  generation of analytic semigroup in $L^2_c$ in the sector $\Sigma_{\pi-\arctan\frac{|a|}{1-|a|}}$ and characterized the domain of the generator through the boundedness in $L^2_c$ of the operators
\begin{align}\label{All-operators} 
		(\lambda- \mathcal{L}_{c,2})^{-1},\quad \Delta_x (\lambda- \mathcal{L}_{c,2})^{-1},\quad \nabla_xD_y (\lambda- \mathcal{L}_{c,2})^{-1}, \qquad \lambda\in \Sigma_{\pi-\arctan\frac{|a|}{1-|a|}}.
\end{align}
On the other hand the multipliers above are bounded in the larger sector  
$$ \Sigma_{\pi-\omega_a}:= \Sigma_{\pi-\arctan\frac{|a|}{\sqrt{1-|a|^2}}}, \qquad \omega_a:=\arctan\frac{|a|}{\sqrt{1-|a|^2}},$$
 by the results in Section \ref{mult}.

In the next lemma we extend the  families  \eqref{All-operators}  to bounded operators on $L^p_m$ on the bigger sector $\Sigma_{\pi-\omega_a}$. In particular  we  prove  that  the resolvents  $(\lambda-\mathcal{L}_{c,2})^{-1}$,  $\lambda\in \Sigma_{\pi-\arctan\frac{|a|}{1-|a|}}$, initially constructed via form method in $L^2_c$, extend to the larger sector $\Sigma_{\pi-\omega_a}$.

We recall the notation of Section \ref{mult}, where  $|a|\leq  \delta\leq 1$,
\begin{align*}
	\theta_\delta=\arctan\frac{|a|}{\sqrt{\delta^2-|a|^2}}\in\left]\arctan\frac{|a|}{\sqrt{1-|a|^2}},\frac\pi 2\right],\qquad \theta_1=\omega_a.
\end{align*}
\begin{lem}\label{lemma multipl}
	Let $0<\frac{m+1}p<c+1$.
	Then the operators  
	$$(\lambda- \mathcal{L}_{c,2})^{-1},\quad \Delta_x (\lambda- \mathcal{L}_{c,2})^{-1},\quad \nabla_xD_y (\lambda- \mathcal{L}_{c,2})^{-1},\quad  B_y (\lambda- \mathcal{L}_{c,2})^{-1}$$ initially defined  on $L^p_m\cap L^2_{c}$ and for $\lambda\in \Sigma_{\pi-\arctan\frac{|a|}{1-|a|}}$, extend for  $\lambda\in \Sigma_{\pi-\omega_a}$ to bounded operators on $L^p_m$,   which we denote   by $\mathcal{R}(\lambda)$, $ \Delta_x \mathcal{R}(\lambda)$, $ \nabla_xD_y \mathcal{R}(\lambda)$, $ B_y \mathcal{R}(\lambda)$, respectively. Moreover 
	for  every $|a|\leq \delta<1$  the family 
	 $\left\{\lambda \mathcal{R}(\lambda):\lambda\in \Sigma_{\pi-\theta_\delta}\right\}$ is $\mathcal{R}$-bounded on $L^p_m$.
\end{lem}
{\sc Proof.}
Let $\lambda\in\Sigma_{\pi-\omega_a}$,    $u \in  C_c^\infty (\R^{N})\otimes  \mathcal{D}$ and $f=\lambda u- \mathcal Lu$.   By taking the Fourier transform with respect to $x$ we obtain
$$
(\lambda+Q_a(\xi)-L_{2a\cdot\xi})\hat u(\xi,\cdot)=\hat f(\xi,\cdot), \qquad  \hat u(\xi, \cdot)=(\lambda+Q_a(\xi)-L_{2a\cdot\xi})^{-1}\hat f(\xi,\cdot).
$$
This means $u={\cal F}^{-1} R_\lambda(\xi) {\cal F} f$, where 
$$R_\lambda(\xi)=(\lambda+Q_a(\xi)-L_{2a\cdot\xi})^{-1}.$$
Theorems \ref{marcinkiewicz} and  \ref{Teo Mult Resolv} yield the  boundedness  of the Fourier multiplier $R_\lambda$ in $L^p\left(\R^N, L^p_m(\R_+)\right)=L^p_m$ and therefore the existence of a bounded operator $\mathcal{R}(\lambda)={\cal F}^{-1} R_\lambda(\xi) {\cal F}\in \mathcal B( L^p_m)$.

Furthermore \cite[Theorem 4.3.9]{Pruss-Simonett} and the $\mathcal R$-boundedness with respect to $\lambda$ of $\lambda R_\lambda(\xi)$ and its $\xi$-derivatives, see again Theorem \ref{Teo Mult Resolv},  imply that  the family the family 
$\left\{\lambda \mathcal{R}(\lambda):\lambda\in \Sigma_{\pi-\theta_\delta}\right\}$ is $\mathcal{R}$-bounded for any $|a|\leq \delta<1$.

Since $ C_c^\infty (\R^{N})\otimes  \mathcal{D}$ is a core for $\mathcal{L}_{c,2}$, this shows in particular that 
\begin{align*}
	(\lambda-\mathcal{L}_{c,2})^{-1}={\cal F}^{-1} R_\lambda(\xi) {\cal F}=\mathcal{R}(\lambda),
\end{align*}
for $\lambda\in \Sigma_{\pi-\arctan\frac{|a|}{1-|a|}}$, where both operators exist. 

However the previous equality although extends to$\lambda\in\Sigma_{\pi-\omega_a}$. In fact the set  $E= \Sigma_{\pi-\omega_a} \cap \rho (\mathcal L_{c,2})$ ($\rho$ denotes the resolvent) is open in  $ \Sigma_{\pi-\omega_a}$ and   $(\lambda-\mathcal{L}_{c,2})^{-1}=\mathcal{R}(\lambda)$ for $\lambda \in E$, by the argument above. But $E$ is also closed in 
 $ \Sigma_{\pi-\omega_a}$ and hence coincides with it. In fact, if $(\lambda_n) \subset E$ converge to $\lambda_0 \in   \Sigma_{\pi-\omega_a}$, then $(\lambda_n-\mathcal L_{c,2})^{-1}=\mathcal{R}(\lambda_n)$ is uniformly bounded and $\lambda_0 \in E$, by elementary spectral theory.

The proof for $\Delta_x \mathcal{R}(\lambda)$, $\nabla_x D_y\mathcal{R}(\lambda)$ are similar. As for  \eqref{Four eq} in  Theorem \ref{Neumnan comp 2}, we have
$$ \Delta_x (\lambda- \mathcal{L}_{c,2})^{-1}=-{\cal F}^{-1} |\xi|^2R_\lambda(\xi) {\cal F}$$  
$$ \nabla_xD_y (\lambda- \mathcal{L}_{c,2})^{-1}=-{\cal F}^{-1} \xi D_yR_\lambda(\xi) {\cal F}$$
and use Theorems \ref{Mult xi^2 Resolv} and 
\ref{Mult D_y/y Resolv} for the  boundedness of the  multipliers  in $L^p_m=L^p(\R^N; L^p_{m}(\R_+))$. 
The boundedness  of $ B_y\mathcal{R}(\lambda)$ follows then by difference.
\qed

\begin{lem} \label{generazione} If  $ 0<\frac{m+1}p<c+1$, an extension $\mathcal L_{m,p}$ of the operator $\mathcal L$, initially defined on $C_c^\infty (\R^{N})\otimes  \mathcal{D}$,  generates an analytic semigroup $\{e^{z\mathcal L_{m,p}}:\ z\in\Sigma_{\frac \pi 2-\omega_a}\}$  in $L^p_m$ which is bounded on  $\Sigma_{\frac\pi 2-\theta_\delta}$, for any $|a|\leq  \delta< 1$. Moreover the semigroup  has maximal regularity and it is consistent with the semigroup generated by $\mathcal L_{c,2}$ in $L^2_{c}$. 
\end{lem}
{\sc Proof.}
Let $|a|\leq  \delta< 1$ and let  us consider the   $\mathcal{R}$-bounded family of operators $\left\{\lambda \mathcal{R}(\lambda):\lambda\in\Sigma_{\pi-\theta_\delta }\right\}$ defined by Lemma \ref{lemma multipl}. In particular it  satisfies
\begin{align*}
	\|\lambda \mathcal{R}(\lambda)\|_{\mathcal{B}(L^p_m)}\leq C,\qquad\forall \lambda\in\Sigma_{\pi-\theta_\delta }.
\end{align*}
By construction  $\mathcal{R}({\lambda})$ coincides with  $(\lambda-\mathcal{L}_{c,2})^{-1}$ when restricted to $L^p_m \cap L^2_{c}$.  Hence, by density,  the family $\left\{\mathcal{R}(\lambda):\lambda\in\Sigma_{\pi-\theta_\delta}\right\}$  satisfies  the resolvent equation
\begin{align*}
	\mathcal{R}(\lambda)-\mathcal{R}(\mu)=(\mu-\lambda)\mathcal{R}(\lambda)\mathcal{R}(\mu),\quad \forall \lambda,\mu\in\Sigma_{\pi-\theta_\delta}
\end{align*}
in $L^p_m$ and therefore it is a pseudoresolvent, see  \cite[Section 4.a]{engel-nagel}. Furthermore $\mbox{rg}(\mathcal{R}(\lambda))$  is dense in $L^p_m$ for every $\lambda \in \Sigma_{\pi-\omega_a}$,  since it contains $ C_c^\infty (\R^{N})\otimes  \mathcal{D}$. 

Let us prove that $\mathcal{R}(\lambda)$ is injective for every $\lambda\in \Sigma_{\pi-\theta_\delta}$. Let $f\in L^p_m$ s.t. $\mathcal{R}(\lambda)f=0$ for some $\lambda \in \Sigma_{\pi-\theta_\delta}$. Since   $\mbox{Ker}(\mathcal{R}(\lambda))=\mbox{Ker}(\mathcal{R}(\mu))$ for any $\lambda,\mu\in\Sigma_{\pi-\theta_\delta}$, see \cite[Lemma 4.5]{engel-nagel},  we have  $\mathcal{R}(\lambda)f=0$ for every $\lambda>0$. Given $\epsilon>0$,  let us choose $g\in L^p_m\cap L^2_{c}$ s.t. $\|f-g\|_{L^p_m}<\epsilon$. Then  
\begin{align*}
	\lambda R(\lambda)g=\lambda R(\lambda)(g-f),\qquad \|\lambda R(\lambda)g\|_{L^p_m}\leq C\epsilon, \qquad \forall\lambda>0.
\end{align*}
Since $\lambda R(\lambda)g=\lambda(\lambda-\mathcal{L}_{c,2})^{-1}g\to g$ as $\lambda\to \infty$ we may suppose, up to  a  subsequence, that $\lambda R(\lambda)g\to g$ a.e.. Then Fatou's Lemma yields
\begin{align*}
	\|g\|_{L^p_m}\leq \liminf_{\lambda\to\infty} \|\lambda \mathcal{R}(\lambda)g\|_{L^p_m}\leq C\epsilon
\end{align*}
which implies  $\|f\|_{L^p_m}\leq \|f-g\|_{L^p_m}+\|g\|_{L^p_m}\leq (1+C)\epsilon$, hence  $f=0$ by the arbitrariness of $\epsilon$  which proves the injectivity of $\mathcal{R}(\lambda)$. 

At this point, the arbitrariness of $\delta$ and   \cite[Proposition 4.6]{engel-nagel} yield the existence of a densely defined closed operator $\mathcal L_{m,p}$  such that $\Sigma_{\pi-\omega_a}\subseteq \rho(\mathcal L_{m,p})$ and $\mathcal R(\lambda)=(\lambda-\mathcal L_{m,p})^{-1}$ for any $\lambda\in\Sigma_{\pi-\omega_a}$. By construction, $(\mathcal L_{m,p};D(\mathcal L_{m,p}))$ extends $\left(\mathcal L, C_c^\infty (\R^{N})\otimes  \mathcal{D}\right)$ and one has
\begin{align*}
	\|\lambda\left(\lambda-\mathcal L_{m,p}\right)^{-1}\|_{\mathcal B\left( L^p_m \right)}\leq C,\qquad \lambda\in \Sigma_{\pi-\theta_\delta},\quad   |a|\leq  \delta< 1.
\end{align*}
Then from standard results on semigroup theory, see for example \cite[Section AII, Theorem 1.14]{nagel}, $(\mathcal L_{m,p},D(\mathcal L_{m,p}))$ generates an analytic semigroup  $\left(e^{z\mathcal{L}_{m,p}}\right)_{z\in \Sigma_{\frac\pi 2-\omega_a}}$   in $L^p_m$ which is bounded on $\Sigma_{\pi-\theta_\delta}$, for any    $|a|\leq  \delta< 1$.

The  maximal regularity of the semigroup follows, using Theorem \ref{MR}, from  the $\mathcal{R}$-boundedness of  the resolvent family $\{\lambda\left(\lambda-\mathcal L_{m,p}\right)^{-1},\  \lambda \in \Sigma_{\pi-\theta_\delta}\}$. Finally, the semigroup is consistent with that in $L^2_{c}$, since the resolvents are consistent.
\qed

We characterize the domain of $\mathcal L_{m,p}$ and collect in one theorem all the results proved in this section.

\begin{teo} \label{generazione1} 
	If $0<\frac{m+1}p<c+1$, then the operator $\mathcal L_{m,p}$  endowed with domain 
	\begin{align*}
		D(\mathcal L_{m,p})&=W^{2,p}_{m,\mathcal{N}}
	\end{align*} 
  generates an analytic semigroup $\{e^{z\mathcal L_{m,p}}:\ z\in\Sigma_{\frac \pi 2-\omega_a}\}$  in $L^p_m$ which is bounded on  $\Sigma_{\frac\pi 2-\theta_\delta}$, for any $|a|\leq  \delta< 1$. Moreover  $C_c^\infty (\R^{N})\otimes \mathcal D$ is a core for $\mathcal L_{m,p}$ and the semigroup  has maximal regularity.
\end{teo}
{\sc Proof.} From Lemma \ref{generazione}  we have only to show that  $D(\mathcal L_{m,p})=W^{2,p}_{m,\mathcal{N}}$. With the notation above, $D(\mathcal L_{m,p})=\mathcal R(1)\left( L^p_m\right )$.  Let $u=\mathcal R(1)f= (I- \mathcal{L}_{c,2})^{-1}f$ with $f \in L^2_{c}\cap L^p_m$. Then 
Lemma \ref{lemma multipl} yields
\begin{align}\label{eq spez}
	\|\Delta_x u\|_{L^p_{m}} +\| B_y u\|_{L^p_{m}}+\|\nabla_xD_y u\|_{L^p_{m}}\leq C\left(\| \mathcal L u\|_{L^p_{m}}+\| u\|_{L^p_{m}}\right).
\end{align}

Using Theorem \ref{Neumnan alpha m 2} and  Theorem \ref{Neumnan comp 2}, we deduce that  $u(x,\cdot)\in D(B^n_{c,2})$ for a.e. $x\in \R^N$. Moreover, $u(x,\cdot)$, $\  B_y u(x, \cdot)\in L^p_m(\R_+)$,  for a.e. $x\in \R^N$. 

Let us show that  $u(x,\cdot)\in D( B^n_{m,p})$. In fact,   setting $f:=u(x,\cdot)-B_y u(x,\cdot)\in L^p_m(\R_+)\cap L^2_{c}(\R_+)$ we have $u=\left(I-B^n\right)^{-1}f\in D(B^n_{m,p})\cap D(B^n_{c,2})$ by the consistency of the resolvent $\left(I-B^n\right)^{-1}$ in $L^p_m(\R_+)$ and  in $L^2_{c}(\R_+)$ .  

Theorem \ref{Neumnan alpha m 2} then implies
\begin{align*}
	\| D_{yy}u\|_{L^p_{m}(\R_+)}+\|y^{-1}D_yu\|_{L^p_{m}}+\|D_{y}u\|_{L^p_{m}(\R_+)}\leq C\|u-  B_y u\|_{L^p_{m}(\R_+)}.
\end{align*}
Then, raising to the power $p$,  integrating over $\R^N$ and using Lemma \ref{lemma multipl} for the last inequality 
\begin{align} \label{eqspez1}
	\|D_{yy}u\|_{L^p_{m}}+\|D_{y}u\|_{L^p_{m}}+\|y^{-1}D_yu\|_{L^p_{m}}
	\leq C \left( \|u\|_{L^p_{m}}+\|B_y u\|_{L^p_{m}}\right )
	\leq C\left( \|u\|_{L^p_{m}}+\|\mathcal {L} u\|_{L^p_{m}}\right ).
\end{align}

By the density of  $L^2_{c}\cap L^p_m$ in $L^p_m$, \eqref{eq spez}, \eqref{eqspez1} hold for every $u \in D(\mathcal L_{m,p})$ and this last is contained in  $W^{2,p}_{m,\mathcal{N}}$, by  \ref{Trace D_yu in W}.

Moreover, since the graph norm is clearly weaker than the norm of $W^{2,p}_{m,\mathcal{N}}$,  \eqref{eq spez}, \eqref{eqspez1} again show that they are equivalent on 
$ D(\mathcal L_{m,p})$, in particular on $C_c^\infty (\R^{N})\otimes \mathcal D$ which is dense in $W^{2,p}_{m,\mathcal{N}}$, by Theorem \ref{core gen}. 

Therefore $ D(\mathcal L_{m,p})=W^{2,p}_{m,\mathcal{N}}$ and in particular $C_c^\infty (\R^{N})\otimes \mathcal D$ is a core.
\qed

\begin{cor} \label{omogeneo}
	Under the hypotheses of Theorem \ref{generazione1} we have for every $u \in W^{2,p}_{m,{\mathcal{N}}}$
	$$
	\| D_{x_i x_j} u\|_{L^p_{m}} +\| D_{x_i y}u\|_{L^p_{m}} +\| D_{yy} u\|_{L^p_{m}}+\|y^{-1} D_{y} u\|_{L^p_{m}}\leq C\| \mathcal Lu\|_{L^p_{m}}.
	$$
\end{cor}
{\sc Proof.} By Theorem \ref{generazione1} the above inequality holds if $\|  u\|_{L^p_{m}}$ is added to the right hand side.
Applying it  to $u_\lambda (x,y)=u(\lambda x, \lambda y)$, $\lambda >0$ we obtain
\begin{align*}
	\| D_{x_i x_j} u\|_{L^p_{m}}+\| D_{x_i y} u\|_{L^p_{m}} +\|D_{yy} u\|_{L^p_{m}}+\|y^{-1} D_{y} u\|_{L^p_{m}}\leq C\left(\| \mathcal L u\|_{L^p_{m}}+\lambda^{-2}\| u\|_{L^p_{m}}\right)
\end{align*}
and the proof follows letting $\lambda \to \infty$.
\qed

\section{General operators and oblique derivative}
Results for more general operators and boundary conditions follow by linear change of variables, as we explain below. Let us first remove the assumption on the special form of  $\mathcal L=\Delta_{x} +2a\cdot\nabla_xD_y+ B_yu$ by considering the general form in  $\R^{N+1}_+$

\begin{equation*}
	\mathcal L =\mbox{Tr }\left(QD^2u\right)+\frac{c}{ y}D_y=\sum_{i,j=1}^{N}q_{ij}D_{x_ix_j}+2\sum_{i=1}^Nq_iD_{x_iy}+\gamma D_{yy}+\frac{c}{y}D_y.
\end{equation*}

If $Q_1$ is the $N\times N$ matrix $(q_{ij})$ and $q=(q_1, \dots, q_N)$ we assume that the quadratic form $Q(\xi,\eta)=Q_1(\xi,\xi)+\gamma \eta^2+2 q\cdot\xi \,\eta$ is positive definite.
Through a linear change of variables in the $x$ variables the term $\sum_{i,j=1}^{N}q_{ij}D_{x_ix_j}$ is transformed into $\gamma\Delta_x$ and all the results of Section 7  hold, replacing $c$ with $ \frac c\gamma$ in the statements (the condition $|a|<1$ of Section 7 is satisfied since the change of variables preserves the ellipticity). The addition of first order terms like $\alpha \cdot \nabla_x+\beta D_y$ is easily treated by standard perturbation theory of analytic semigroups and maximal regularity; the case of variable coefficients can also be handled by freezing the coefficients and will be done in the future to deal with degenerate problems in bounded domains.


A further change of variables allows to deal with  the operator 
\begin{equation*}
	\mathcal L =\mbox{Tr }\left(QD^2u\right)+\frac{ v\cdot \nabla }{y}=\sum_{i,j=1}^{N+1}q_{ij}D_{ij}+\frac{b\cdot \nabla_x}{y}+\frac{cD_y}{y},\qquad c\neq 0, 
\end{equation*}
where $v=(b,c)$ and $Q$ positive definite. We impose an oblique derivative boundary condition  $v \cdot \nabla u(x,0)=0$ in the integral form
\begin{align*}
\frac{v\cdot \nabla u}y=\frac{b\cdot \nabla_x u+c D_yu}y\in L^p_m
\end{align*}
and define therefore
\begin{align*}
	W^{2,p}_{m,v}=\{u \in W^{2,p}_{m}:\  y^{-1}v\cdot \nabla u \in L^p_m\}.
\end{align*}
We transform $\mathcal L$ into a similar operator with $b=0$ and Neumann boundary conditions by defining the following isometry of $L^p_m$
\begin{align}\label{Tran def}
	T\, u(x,y)&:=u\left(x-\frac b c y,y\right),\quad (x,y)\in\R_+^{N+1}.
\end{align}

\begin{lem}\label{Isometry action der} Let $1< p< \infty$,  $ v=(b,c)\in\R^{N+1}$,  $c\neq 0$. Then  for $u\in W^{2,1}_{loc}\left(\R_+^{N+1}\right)$ 
	\begin{itemize}
	\item[(i)] 
	\begin{align*}
		 T^{-1}\,\left(\mbox{Tr }\left(QD^2u\right)+\frac{ v\cdot \nabla u }{y}\right) T\,u=\mbox{Tr }\left(\tilde QD^2u\right)+\frac{c}{y}D_y u
	\end{align*}
where $\tilde Q$  is a uniformly elliptic symmetric matrix defined by 
\begin{align*}
	\tilde Q=
	\left(
	\begin{array}{c|c}
		Q_N-\frac{2}c b\otimes q-\frac \gamma{c^2} b\otimes b  & q^t- \frac{\gamma}c b^t \\[1ex] \hline
		q- \frac{\gamma}c b & \gamma
	\end{array}\right)
\end{align*} 
and  $\gamma= q_{N+1,N+1}$.
	\item[(ii)] $\ds T\,\left(W^{2,p}_{m, \mathcal N}\right)	=W^{2,p}_{m,  v}$.
	\end{itemize}
\end{lem}{\sc{Proof.}} The proof follows by a straightforward computation.
\qed

We can therefore deduce  results also for the last operator which we state only in the parabolic setting. The proof follows directly from the above lemma, the general theory  of Section 2 and standard semigroup theory.

\begin{teo}
Let $0<\frac{m+1}{p}<\frac{c}\gamma +1$, $v=(b,c)$ with $c \neq 0$, $Q$ uniformly elliptic and
\begin{equation*}
	\mathcal L =\mbox{Tr }\left(QD^2u\right)+\frac{ v\cdot \nabla }{y}
\end{equation*}
with domain  $W^{2,p}_{m,v}$. Then for each $1<q<\infty$, $T>0$ and $u_0 \in W^{2,p}_{m, v}, f\in L^q([0,T];L^p_m)$ the problem
$$\frac{\partial}{\partial t} u(t,x,y)-\mathcal Lu(t,x,y)=f(t,x,y), \quad t>0, \qquad u(0,x,y)=u_0(x,y)$$
admits a unique solution  $u\in W^{1,q}([0,T];L^p_m)\cap L^q([0,T];W^{2,p}_{m,v})$.

\end{teo}
\bibliography{../../TexBibliografiaUnica/References}

\end{document}